\documentclass[11pt,a4paper]{article}

\usepackage[margin=2.5cm]{geometry}
\usepackage{amsmath,amssymb,amsthm}
\usepackage{enumitem}
\usepackage{booktabs}
\usepackage{graphicx}
\usepackage{hyperref}
\hypersetup{colorlinks=true,linkcolor=blue!55!black,urlcolor=blue!55!black,
            citecolor=blue!55!black}
\usepackage{xcolor}
\usepackage{tikz}
\usetikzlibrary{calc,decorations.pathmorphing}
\definecolor{batA}{HTML}{C8102E}   
\definecolor{batB}{HTML}{0098D4}   
\definecolor{batC}{HTML}{0B2A8C}   
\definecolor{batD}{HTML}{00A040}   
\tikzset{
  bat1/.style={draw=batA,line width=2.2pt},
  bat2/.style={draw=batB,line width=1.8pt,line cap=round,dash pattern=on 0pt off 2.7pt},
  bat3/.style={draw=batC,line width=1.3pt},
  bat4/.style={draw=batD,line width=2.2pt,line cap=round,dash pattern=on 0pt off 2.65pt},
  photon/.style={decorate,decoration={snake,amplitude=2.2pt,segment length=9pt,
    pre length=3pt,post length=3pt}},
  spring/.style={decorate,decoration={coil,aspect=0.5,amplitude=4.5pt,
    segment length=8pt,pre length=3pt,post length=3pt}}}
\DeclareRobustCommand{\photonsample}{\tikz[baseline=-0.6ex]
  \draw[line width=1.2pt,photon] (0,0) -- (1.1,0);}
\DeclareRobustCommand{\springsample}{\tikz[baseline=-0.6ex]
  \draw[line width=1.2pt,spring] (0,0) -- (1.1,0);}

\theoremstyle{plain}
\newtheorem{theorem}{Theorem}[section]

\newtheorem{lemma}[theorem]{Lemma}
\newtheorem{corollary}[theorem]{Corollary}
\theoremstyle{definition}
\newtheorem{definition}[theorem]{Definition}
\newtheorem{remark}[theorem]{Remark}

\newcommand{\R}{\mathbb{R}}
\newcommand{\Z}{\mathbb{Z}}
\newcommand{\ms}[1]{\{\!\{\,#1\,\}\!\}}
\newcommand{\norm}[1]{\lVert #1\rVert}
\newcommand{\pt}[1]{[#1]}
\DeclareMathOperator{\sort}{sort}

\newcommand{\MPS}{Matou\v{s}ek, P\v{r}\'{\i}v\v{e}tiv\'{y} and \v{S}kovro\v{n}}

\title{How many sorted projections determine any $n$ points in $\R^k$?}
\author{%
  Christopher Gorham Lester\thanks{Corresponding author.}\\
  {\small\href{mailto:CGL20@cam.ac.uk}{\texttt{CGL20@cam.ac.uk}}}
  \and
  Raymond Hau S{\o}rensen\\
  {\small\href{mailto:RS2398@cam.ac.uk}{\texttt{RS2398@cam.ac.uk}}}
  \and
  {\small Cavendish Laboratory, Department of Physics, University of Cambridge, UK}}
\date{12 November 2026}

\begin{document}
\maketitle

\begin{abstract}
\noindent
Let $d_1,\dots,d_m$ be directions in $\R^k$, where $k\ge2$.  A multiset $X$ of $n$ points of
$\R^k$ projects to the $m$ multisets of numbers $\ms{x\cdot d_j : x\in X}$,
or equivalently to the $m$ sorted lists of these numbers.  How many directions
are needed for these projections to determine every multiset of $n$ points?
Two simple bounds are known in every dimension: more than $(k-1)\log_2n$
directions are always needed, and $(k-1)n+1$ directions in general position
always suffice.  In the plane much more is known: \MPS{} showed that almost
every choice of $O(\log n\,\log\log n)$ directions suffices.  Whether anything
similar holds in higher dimensions was left open by Dym, Wellershoff,
Tsoukanis, Levy and Balan.  We observe that two lemmas of Hair and Sahai, proved for list
recovery of random linear codes, answer that question, and by tightening their argument
we show that in every dimension almost every choice of
$(4k-2)\ln(2n)\ln(k\log_2 2n)+k+1$ directions suffices.  For $k\le10$ this beats
the linear bound once $n$ is larger than about $100$: for example, $144$
directions suffice for multisets of $1000$ points in the plane, and an exact
computation lowers this to $133$.  Letting some of the random blocks in their
argument fail removes a factor $\log k$: almost every choice of
$O(k\log n\,\log\log n)$ directions suffices, uniformly in $k$ and $n$.  We also show that $2k-1$ directions never determine three points,
so that the logarithmic lower bound is not tight; that every $m$ directions in
the plane fail on some pair of multisets of size $O(1.65^m)$, improving the
bound $O(1.80^m)$ of \MPS; and that for each $k$, $m$ and $n$, either almost
every choice of $m$ directions works or almost none does.
\end{abstract}

\section{Introduction}
\label{sec:intro}

\subsection{The problem}

For a nonzero vector $d\in\R^k$ and a finite multiset $X$ of points of $\R^k$,
the \emph{projection} of $X$ on $d$ is the multiset of real numbers
\[
  \pi_dX\;:=\;\ms{\,x\cdot d \;:\; x\in X\,}.
\]
For a list $D=(d_1,\dots,d_m)$ of nonzero vectors, the \emph{directions}, let
\[
  \Phi_D(X)\;:=\;\bigl(\sort\pi_{d_1}X,\,\dots,\,\sort\pi_{d_m}X\bigr)
  \;\in\;\R^{mn},
\]
where $n=|X|$ and $\sort$ lists a multiset of $n$ numbers in increasing order.
Only the line through a direction matters, since $\pi_{\lambda d}X$ is $\pi_dX$ scaled by $\lambda$.

We say that $D$ \emph{separates} multisets of size $n$ if $\Phi_D$ is injective
on them, that is, if any two multisets of $n$ points with the same projection
on every $d_j$ are equal.  We want to know two things: how many directions are needed to achieve this injectivity, and how should those directions be chosen.  We denote the answer to the first question by $m_k(n)$, i.e.:
\[
  m_k(n)\;:=\;\min\{\,m : \text{some list of $m$ directions in $\R^k$
  separates multisets of size $n$}\,\}
\]
and we assume $k\ge 2$ throughout this paper, as $k=1$ is trivial.\footnote{There is only one direction in $k=1$ dimensions, and it separates multisets of every size, so $m_1(n)=1$ for all $n$.}  
The answers to these questions are of interest in two areas that interest us:

\emph{Sorting-based embeddings.}  In machine learning, $\Phi_D$ is used to
present an unordered set of vectors to an algorithm that expects a vector of
fixed length; it is the map $\beta_A$ of \cite{BHS,DWTLB}, the directions being
the columns of the matrix $A$.  It does not depend on the order in which the
points are listed, and it is piecewise linear.  When it is injective it is
bi-Lipschitz \cite[Theorem~1.2]{BHS}, whereas no invariant map that is
differentiable at a multiset with a repeated point can be \cite[Theorem~21]{CIM}.  In particular, the sum-pooling embedding
$X\mapsto\sum_{x\in X}f(x)$ of Deep Sets \cite{ZKRPSS} is never bi-Lipschitz if
$f$ is differentiable at even one point \cite[Proposition~5.1]{AGARD}.  It costs
about $mn(k+\log_2n)$ operations to evaluate, so few directions mean a fast
embedding.

\emph{Discrete tomography.}  In the plane, $\pi_dX$ records how many points of
$X$ lie on each line orthogonal to $d$.  This is the \emph{discrete X-ray} of
$X$ in the direction $d^\perp$, and $m_2(n)$ is the question of \MPS{}
\cite{MPS}, how many X-rays determine every set of $n$ points, asked for
multisets, which they also treat.

The difficulty is the order of the quantifiers.  For a fixed multiset $X$,
$k+1$ random directions determine $X$ with probability one if they are
independent and each is orthogonal to any given nonzero vector with probability
zero, as when each is uniform on the unit sphere \cite{TFD}; and a suitable
fixed list of $k+1$ directions determines almost every multiset
\cite[Proposition~3.7]{BHS}.  But probability one is not enough!  We ask
for one list of directions that works for every multiset of size $n$ in $k$-dimensions at once. Failure on a measure zero set of points is not permissible if you want injectivity, and this is what makes the problem much harder.

\begin{figure}[tp]
\centering
\includegraphics[width=\textwidth]{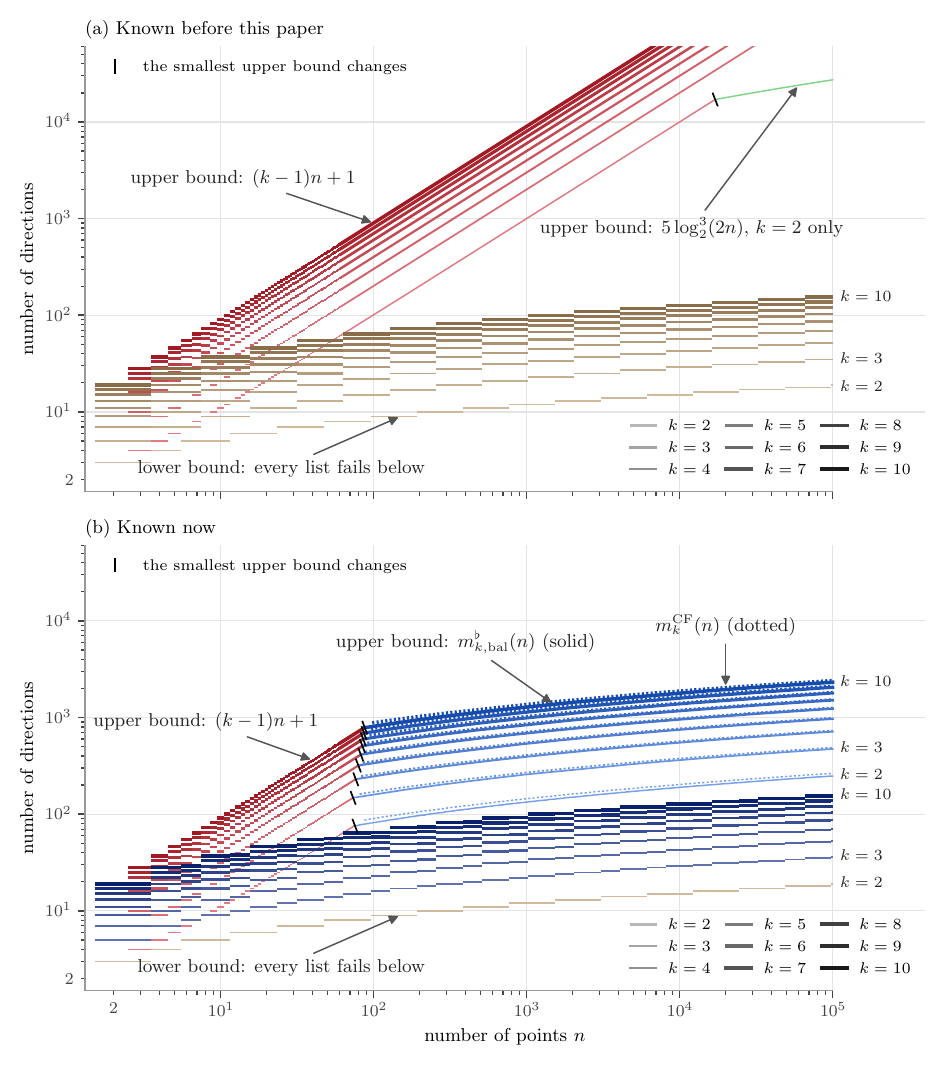}
\caption{Bounds on the number of directions that separate multisets of $n$
points in $\R^k$, for $k=2,\dots,10$ (thin and light to thick and dark):
(a)~before this paper and (b)~now; this paper's bounds are blue.  Every list
of fewer directions than the lower curve fails \cite[Theorem~5]{DWTLB} (for
$k=2$, \cite{BL90} and \cite{MPS}), in~(b) for $k\ge3$ also by
Theorem~\ref{thm:three}.  The upper curve is the smallest known upper
bound: the linear bound $(k-1)n+1$, for every list in general position
\cite{Renyi,HSW}; for almost every list, in~(a) for $k=2$ the bound
$5\log_2^3(2n)$ of \cite[Proposition~3.1]{MPS}, and in~(b)
$m^\flat_{k,\mathrm{bal}}(n)$ of Theorem~\ref{thm:main}, with the closed form
$m^{\mathrm{CF}}_k(n)$ dotted (on this range the smaller of $m^{\mathrm{CF}}_k(n)$
and $m^{\mathrm{HD}}_k(n)$, Remark~\ref{rem:no-logk}).
Better lower bounds for $k=2$ (\cite{MPS} and Theorem~\ref{thm:plane}) and the solved cases with small $n$ are
in Figures~\ref{fig:best-k-before} and~\ref{fig:best-k-now}, and
Table~\ref{tab:numbers} gives numbers.  Drawn by \texttt{plot\_bounds.py} in
the accompanying code.}
\label{fig:best}
\end{figure}

\begin{figure}[p]
\centering
\includegraphics[width=\textwidth]{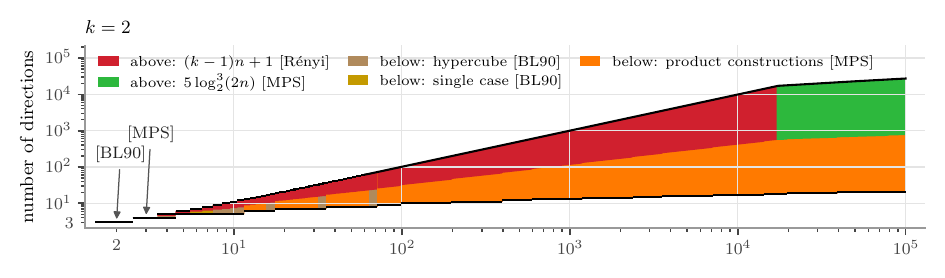}\\
\includegraphics[width=\textwidth]{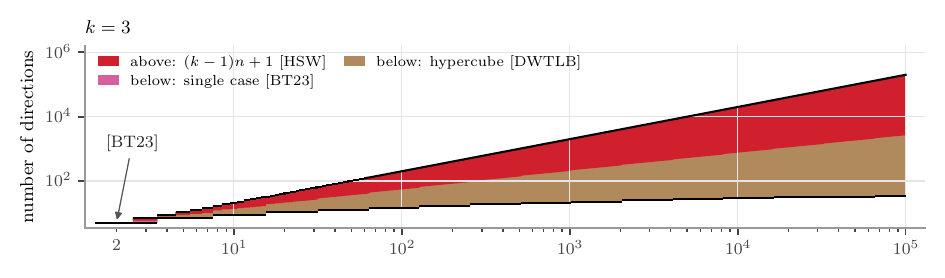}\\
\includegraphics[width=\textwidth]{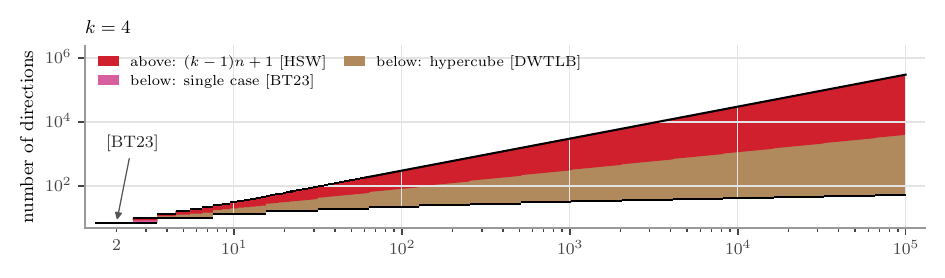}\\
\includegraphics[width=\textwidth]{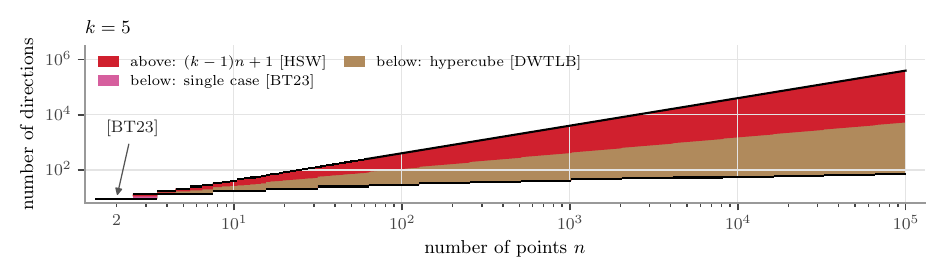}
\caption{The best bounds on $m_k(n)$ known before this paper, for
$k=2,\dots,5$, including those for single values of $n$.  $m_k(n)$ lies in the
shaded range between the upper and lower bounds (black).  Each half of the
shading has the colour of the bound credited on its side, named in the legend,
and where the two meet a label names the reference credited with $m_k(n)$.  The
credits follow fixed rules from a database of what each reference states
(\texttt{bounds.toml} in the accompanying code): a bound goes to the earliest
reference that states it.  The figures show only results that come with a proof
a reader can check; findings made in floating point are acknowledged in the
text instead (Remark~\ref{rem:small} and the provenance note).}
\label{fig:best-k-before}
\end{figure}

\begin{figure}[p]
\centering
\includegraphics[width=\textwidth]{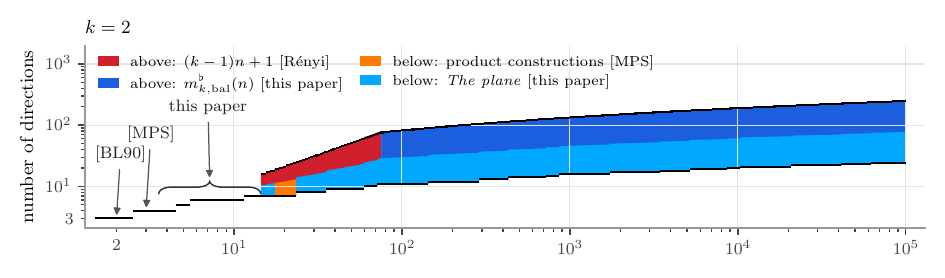}\\
\includegraphics[width=\textwidth]{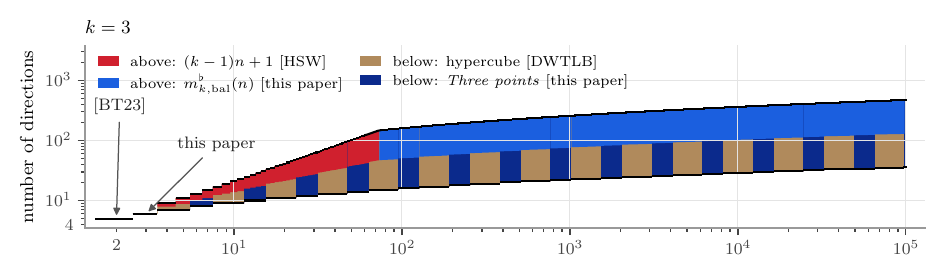}\\
\includegraphics[width=\textwidth]{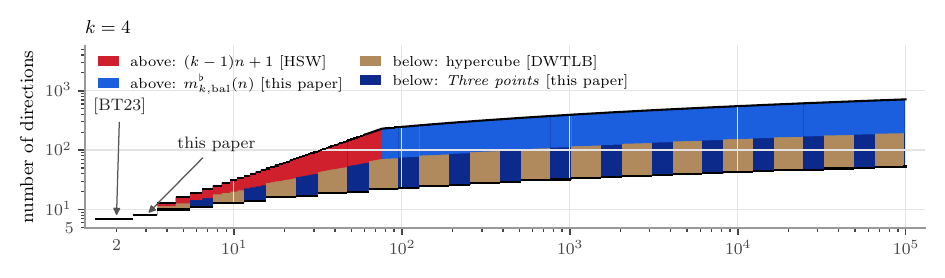}\\
\includegraphics[width=\textwidth]{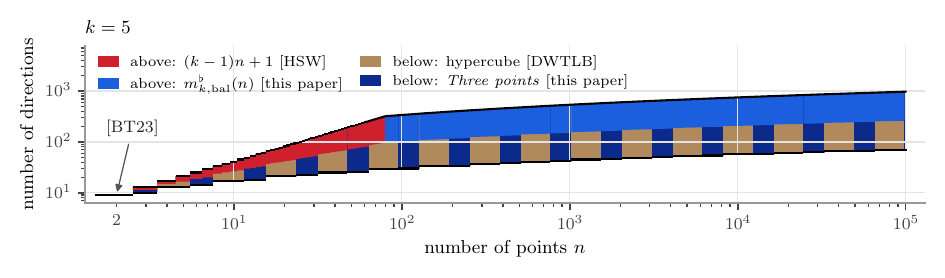}
\caption{As Figure~\ref{fig:best-k-before}, with the results of this
paper.  \emph{Three points} is Theorem~\ref{thm:three}, \emph{The plane} is
Theorem~\ref{thm:plane}, and $m^\flat_{k,\mathrm{bal}}(n)$ is from
Theorem~\ref{thm:main}.
Values that others found in floating point and we verified exactly
(Remark~\ref{rem:small}) are credited to this paper.  For $k=2$ the solved
cases with $6\le n\le14$ are from Remark~\ref{rem:six-seven}, and for
$18\le n\le23$ the lower bound credited to the product constructions of
\cite{MPS} has the same values as \emph{The plane}, which it predates.}
\label{fig:best-k-now}
\end{figure}

\subsection{What is already known}

\begin{itemize}[leftmargin=*,itemsep=2pt]
\item \emph{A logarithmic lower bound.}  Every list of $m$ directions fails
  when $\lceil m/(k-1)\rceil\le\log_2n+1$ \cite[Theorem~5]{DWTLB}; for $k=2$
  it is Proposition~4 of Bianchi and Longinetti \cite{BL90}, and \cite{MPS}
  note that the construction has been rediscovered several times.  So
  $m_k(n)\ge(k-1)\bigl(\lfloor\log_2n\rfloor+1\bigr)+1$.  For $n=2$ this is
  $2k-1$, which the linear bound attains; Balan and Tsoukanis \cite{BT23} had
  shown that $m_k(2)=2k-1$, by reduction to real phase retrieval.
\item \emph{A linear upper bound.}  Every list of $(k-1)n+1$ directions in
  general position (any $k$ of them linearly independent,
  Definition~\ref{def:gp}) separates multisets of size $n$.  For $k=2$ this is due to R\'enyi \cite{Renyi}.  In every
  dimension it was proved by Holler, Schl\"uter and Wirth
  \cite[Proposition~48 and Corollary~49]{HSW}, and found again in
  \cite[Theorem~4]{DWTLB}.
\item \emph{The plane.}  \MPS{} \cite{MPS} showed that almost every list of
  $m$ directions in the plane separates multisets of size up to
  $2^{cm/\log m}$, for an absolute constant $c>0$.  So
  $m_2(n)=O(\log n\,\log\log n)$, with constants that they do not make
  explicit.  A weaker form of their argument is explicit
  \cite[Proposition~3.1 and Lemma~2.3]{MPS}: almost every list of
  $5\log_2^3(2n)$ directions in the plane separates multisets of size $n$, which
  beats the linear bound from $n=17\,078$.  In the other direction they constructed, for every list of $m$
  directions, two different multisets of size $O(1.7996^m)$ with the same
  projections.
\item \emph{Codes.}  In a paper on list recovery of random linear codes over
  finite fields, Hair and Sahai \cite{HS26} carried the method of \cite{MPS}
  to $d$ dimensions.  They state only orders of growth, but two of their
  lemmas are exactly what our problem needs.
\end{itemize}
Dym et al.\ asked whether the planar bound extends to higher dimensions:
``it remains unclear whether similar bounds hold when $d>2$''
\cite[Remark~6]{DWTLB} (their $d$ is our $k$).  The lemmas of \cite{HS26}
settle this, as we explain below.

\subsection{What we will show}

We prove an upper bound on the number of directions needed
(Theorem~\ref{thm:main}), two lower bounds (Theorems~\ref{thm:three}
and~\ref{thm:plane}), and a zero--one law that compares almost every list of
directions with particular lists (Theorem~\ref{thm:zero-one}).

\paragraph{Upper bounds.}  We say that a property holds for \emph{almost every} list of $m$ directions if
the lists that fail it form a set of Lebesgue measure zero in $\R^{mk}$, so
that directions drawn at random with a density, on $\R^k$ or on the unit
sphere, have it with probability one.  Appendix~\ref{app:conditions} defines
this and the other conditions that we put on the directions, says which result
uses which, and compares ways of drawing directions at random.

\begin{theorem}[Few directions suffice]
\label{thm:main}
Let $k\ge2$ and $n\ge2$.  For $m\ge2$ and $1\le t\le\lceil m/2\rceil$, write
\[
  P_{2n}(m,t)\;:=\;\log_2(2n)\prod_{j=0}^{t-1}\Bigl(1-\frac{m}{4(m-j)\ln(2n)}\Bigr),
\]
and for $m\ge q\ge k$, write $m=qt+r$ with $0\le r<q$ and
$B_{2n}(m,q):=(q-r)P_{2n}(m,t)+rP_{2n}(m,t+1)$.  Define
\begin{subequations}
\begin{align}
\label{eq:msharp}
  m^\sharp_k(n)&\;:=\;k\cdot\min\Bigl\{T\ge1:\
  k\,P_{2n}(kT,T)<1\Bigr\},\\[2pt]
\label{eq:mflateq}
  m^\flat_{k,\mathrm{eq}}(n)&\;:=\;\min\Bigl\{qT:\ q\ge k,\ T\ge1,\
  q\,P_{2n}(qT,T)<q-k+1\Bigr\},\\[2pt]
\label{eq:mflatbal}
  m^\flat_{k,\mathrm{bal}}(n)&\;:=\;\min\Bigl\{m\ge k:\
  B_{2n}(m,q)<q-k+1 \text{ for some integer } q\in\{k,\dots,m\}\Bigr\}.
\end{align}
\end{subequations}
Then
\[
  m^\flat_{k,\mathrm{bal}}(n)\;\le\;m^\flat_{k,\mathrm{eq}}(n)\;\le\;m^\sharp_k(n)\;\le\;\bigl\lceil m^{\mathrm{CF}}_k(n)\bigr\rceil,
\]
where
\begin{equation}
\label{eq:mcf}
  m^{\mathrm{CF}}_k(n)\;:=\;(4k-2)\ln(2n)\,\ln\bigl(k\log_2(2n)\bigr)+k+1 .
\end{equation}
If $m\ge m^\flat_{k,\mathrm{bal}}(n)$ (and in particular if $m$ is at least
$m^\flat_{k,\mathrm{eq}}(n)$, $m^\sharp_k(n)$, or $m^{\mathrm{CF}}_k(n)$), then
almost every list (Definition~\ref{def:ae}) of $m$ directions in $\R^k$
separates multisets of size $n$.  Hence, by Remark~\ref{rem:no-logk},
$m_k(n)=O(k\log n\,\log\log n)$.
\end{theorem}

The proof tightens an argument of Hair and Sahai \cite[Lemmas~5--7]{HS26}, who
carried the random-colour method of \cite{MPS} to $d$ dimensions for list
recovery of random linear codes.  Their lemmas apply to sorted projections with
almost no work, and as their proof is written they give the bound
$m^{\mathrm{HS}}_k(n):=k\bigl\lceil8\ln(2n)\bigl(\ln\ln(8n)+\ln(2k)+1\bigr)
\bigr\rceil$ (Section~\ref{sec:graphs}).  That observation, which brings their
work to sorting embeddings, is ours, but is built entirely on their work.  For $k=2$ and
$n=1000$, for example, $m^{\mathrm{HS}}_2(n)=558$, and their own final inequality,
evaluated exactly with nothing else changed, gives $m^{\mathrm{HSE}}_2(n)=458$
(equation~\eqref{eq:mhse}).  Their proof gives away
constant factors in four places, and keeping them (Lemma~\ref{lem:sharp}) is
our contribution to the upper bound: it gives $m^\sharp_2(n)=134$.  Dividing
the colours into $q\ge k$ blocks instead of $k$ and asking only that at least
$k$ of the $q$ blocks succeed (Lemma~\ref{lem:blocks-fail}) gives the
\emph{equal-block bound} $m^\flat_{k,\mathrm{eq}}(n)$ when the blocks have
equal size $T$, and the \emph{balanced-block bound}
$m^\flat_{k,\mathrm{bal}}(n)$ when the $m$ colours are shared among the $q$
blocks as evenly as possible ($\lfloor m/q\rfloor$ or $\lceil m/q\rceil$ per
block), lowering $m^\sharp_2(1000)=134$ to $m^\flat_{2,\mathrm{bal}}(1000)=133$,
with larger savings as $k$ grows (for $n=100$ and $k=20$,
$m^\sharp_{20}(100)=2040$, $m^\flat_{20,\mathrm{eq}}(100)=1752$, and
$m^\flat_{20,\mathrm{bal}}(100)=1751$; Table~\ref{tab:flat-comparison}).
Taking $q=2k$ in $m^\flat_{k,\mathrm{eq}}(n)$ also gives a simple
high-dimensional bound
$m^{\mathrm{HD}}_k(n):=2k\bigl\lceil 4\ln(2n)\ln\bigl(2\log_2(2n)\bigr)\bigr\rceil$
that removes the $\log k$ term and scales as $O(k\log n\,\log\log n)$ when $k$
is large relative to $n$ (Remark~\ref{rem:no-logk}).  For fixed $k$ the theorem
matches the planar result of \cite{MPS}, and so answers the question of
\cite{DWTLB}.

For small $n$ the linear bound $(k-1)n+1$ is smaller, and it holds for every
list in general position, not only for almost every list.  The theorem improves
on it for all $n\ge76$ when $k=2$, $n\ge75$ ($74$ for $m^\flat_{k,\mathrm{bal}}$)
when $k=3$ and $n\ge94$ ($88$ for $m^\flat_{k,\mathrm{eq}}$ and
$m^\flat_{k,\mathrm{bal}}$) when $k=10$; beyond that, $m^\flat_{k,\mathrm{bal}}(n)$
is the smallest upper bound we know, and a short computation gives it exactly.
The bound $m^{\mathrm{CF}}_k(n)$ of the theorem
is the smallest we know in closed form for moderate $k$; for a quick bound
without iteration, use $\min\{m^{\mathrm{CF}}_k(n),m^{\mathrm{HD}}_k(n)\}$
(Remark~\ref{rem:no-logk}), which is $m^{\mathrm{CF}}_k(n)$ unless $k$ exceeds
about $4\log_2(2n)$.  It is only slightly larger, $144$
instead of $134$ for $k=2$ and $n=1000$, and it improves on the linear bound
from $n=87$, $81$ and $98$ respectively.  Figure~\ref{fig:best} shows the best bounds known
before this paper and now, Figures~\ref{fig:best-k-before}
and~\ref{fig:best-k-now} add the solved small cases for $k=2,\dots,5$, and
Tables~\ref{tab:numbers} and~\ref{tab:flat-comparison}, at the end of the paper,
give the bounds as numbers.

Theorem~\ref{thm:main} fixes $n$ and $m$ before choosing the directions.  In an
application the size of the multisets may not be known in advance, or may grow,
and one would like to draw the directions once and simply use more of them for
larger multisets, without redrawing those already in use.  This is possible:
with probability one, a single random sequence $d_1,d_2,\dots$ has the property
that, for every $n$, its first $m$ directions separate multisets of size $n$
whenever $m\ge m^\flat_{k,\mathrm{bal}}(n)$ (and in particular whenever
$m\ge m^\sharp_k(n)$; Corollary~\ref{cor:nested}).  Nor does it hurt to fix the
first $k$ directions to be the coordinate axes and draw only the others
(Corollary~\ref{cor:axes}).

\paragraph{Lower bounds.}  For these it is convenient to count from the other
side.  For a list
$D$ of directions let $f(D)$ be the smallest size of two different multisets
with the same projection on every direction of $D$ (or $\infty$ if there are
none).  Then $D$ separates multisets of size $n$ exactly when $n<f(D)$
(Lemma~\ref{lem:pad}).

\begin{theorem}[Three points]
\label{thm:three}
Let $k\ge2$.  Every list $D$ of $2k-1$ directions in $\R^k$ has $f(D)\le3$:
it does not separate multisets of size three.  Consequently:
\begin{enumerate}[label=(\alph*)]
\item $2k\le m_k(3)\le3k-2$;
\item if $3\cdot2^{q-2}\le n<2^q$ for some integer $q\ge2$, then
  $m_k(n)\ge(k-1)\bigl(\lfloor\log_2n\rfloor+1\bigr)+2$, one more than the
  logarithmic lower bound.\footnote{The logarithmic lower bound, or one more
    than it when Theorem~\ref{thm:three}(b) applies, may be written jointly as
    $m_k(n)\ge1+(k-2)\,\ell(n)+\ell(4n/3)$, where $\ell(x)=\lfloor\log_2x\rfloor+1$.}
\end{enumerate}
\end{theorem}

For $k=2$ the first statement is in \cite{MPS}.  Exact certificates
(Remark~\ref{rem:small}) give $m_3(3)=6$ and $m_4(3)=8$, so the lower bound
in~(a) is attained there.  Dym et al.\ \cite[Section~4]{DWTLB} had found such
a list for $m_4(3)$ in floating point; the one they print for $m_3(3)$ does not
separate.

\begin{theorem}[The plane]
\label{thm:plane}
Write $m=5t+r$ with $0\le r\le4$.  Every list $D$ of $m\ge1$ directions in
$\R^2$ has $f(D)\le\frac12\,c_r\,12^t$, where $(c_0,\dots,c_4)=(1,2,4,6,10)$.
Hence no list of $m$ directions in the plane separates multisets of size
$\frac12\,12^{\lceil m/5\rceil}=O(1.6438^m)$, and
$m_2(n)>5\log_{12}(2n)-5\approx1.395\log_2n-3.6$.
\end{theorem}

This improves the bound $O(1.7996^m)$ of \cite{MPS}.  It combines their product
construction with a six-point example for every list of five directions
(Lemma~\ref{lem:five}) and, when $r=4$, with a five-point example of Bianchi
and Longinetti \cite{BL90} for every list of four.\footnote{We are grateful to
\cite{ZM26} for suggesting, after reading the first public version of this
paper, that the result we gave at that time (that any four directions always
fail on five points in the plane) could be improved by one on each count.  This
reminded us to follow up our remark ``5 generic directions are fooled by a 6+6
pair; 6 are not'' (made in \cite{LS26}) with our own proof of the stronger
proposition.  No
proof was provided by \cite{ZM26}, but there is no reason to believe they did
not have one.}  The nine-direction example of \cite{MPS} ($f\le99$) never
helps.  Nor do six directions, by an exhaustive computer search that also
settles $m_2(n)$ for every $n\le14$ (Remark~\ref{rem:six-seven}).

\paragraph{Almost every list behaves alike.}  Theorem~\ref{thm:main}
speaks of almost every list.  The next theorem shows that, for given $k$, $m$
and $n$, almost every list behaves in the same way, and says which particular
lists are sure to behave like them.

\begin{theorem}[A zero--one law]
\label{thm:zero-one}
Fix $k$, $m$ and $n$.  Either almost every list of $m$ directions in $\R^k$
separates multisets of size $n$, or almost every list does not.  Every list
whose $mk$ coordinates are algebraically independent (Definition~\ref{def:ai})
behaves like almost every list.
\end{theorem}

In the plane, \MPS{} note that their bound $2^{cm/\log m}$ holds for every list
with algebraically independent coordinates \cite[Section~5, footnote~4]{MPS}.  By the first
statement, a list that separates with fewer directions than almost every list
needs, if one exists, cannot be found by drawing directions at random with a
density (Table~\ref{tab:random}).  We do
not know whether such lists exist.

\subsection{Method}

The proof of Theorem~\ref{thm:main} follows the strategy of \cite{MPS}, which
Hair and Sahai \cite{HS26} carried to $d$ dimensions in another setting.
Suppose that $D$ fails, so that there are two different multisets $X$ and $Y$
of size at most $n$, with no point in common, that have the same projection on
every direction.  For each direction $d_j$, pairing each point of $X$ with a
point of $Y$ that has the same projection on $d_j$ gives a perfect matching
between $X$ and $Y$.  The union of these $m$ matchings, with the edges of the
$j$-th matching given colour $j$, is the \emph{matching graph} of the pair.
Each edge $\{a,b\}$ of colour $j$ is a linear equation $(a-b)\cdot d_j=0$
satisfied by the points.
\begin{itemize}[leftmargin=*,itemsep=2pt]
\item \emph{Rigidity} (Section~\ref{sec:rigidity}).  Suppose that, on some set
  $W$ of at least two vertices, the matching graph has $k$ spanning trees that
  use pairwise disjoint sets of colours.  Then the coordinates of the
  directions satisfy a polynomial equation that is not identically zero
  \cite[Lemma~7]{HS26}.  So for almost every $D$ this never happens.
\item \emph{Graph theory} (Section~\ref{sec:graphs}).  Every edge-coloured
  multigraph on $N$ vertices whose colour classes are perfect matchings
  contains such a pattern once it has enough colours.  The randomised argument
  of \cite[Lemma~6]{HS26} needs about $8k\ln N\,(\ln\ln N+\ln k)$ colours as
  written; with the constants that it gives away kept, about
  $(4k-2)\ln N\,\ln(k\log_2N)$ (Lemma~\ref{lem:sharp}).
\end{itemize}
For $k=2$ both steps are in \cite{MPS}.

The lower bounds are constructions.  If $D_1$ and $D_2$ fail on pairs of sizes
$s_1$ and $s_2$, then the concatenated list fails on a pair of size $2s_1s_2$
\cite[Lemma~4.1]{MPS}.  Theorem~\ref{thm:three} comes from a cancellation in
the resulting ``hypercube'' construction, and Theorem~\ref{thm:plane} from
examples for every list of five and of four directions (Lemmas~\ref{lem:five}
and~\ref{lem:four}).  Theorem~\ref{thm:zero-one} is a standard argument
applied to a criterion for separation from \cite{BHS}
(Lemma~\ref{lem:criterion}).

\paragraph{Organisation.}  Section~\ref{sec:prelim} collects the facts we
use about pairs with equal projections.  Section~\ref{sec:lower} proves
lower bounds, Theorems~\ref{thm:three} and~\ref{thm:plane};
Sections~\ref{sec:rigidity} and~\ref{sec:graphs} prove upper bounds,
Theorem~\ref{thm:main}; and Section~\ref{sec:generic} proves
Theorem~\ref{thm:zero-one} and its consequences.  Section~\ref{sec:open} lists open questions, and
Appendix~\ref{app:conditions} collects the conditions that the results put on the directions, and Appendix~\ref{app:exclusion} draws the results for $k=2$ and $k=3$ as exclusion plots.  The computations are in the accompanying code, which is in the ancillary files of this paper's arXiv submission.

\section{Preliminaries}
\label{sec:prelim}

Fix $k\ge2$ and a list $D$ of directions.  Two different multisets $X$, $Y$
of the same size $s$ with $\pi_dX=\pi_dY$ for every $d$ in $D$ form a
\emph{confusable pair} of size $s$, and $f(D)$ is the smallest size of a
confusable pair.  We write $X\uplus Z$ for the sum of multisets
(multiplicities add) and $X\cap Y$ for their common part.

\begin{lemma}[Padding and cancelling]
\label{lem:pad}
If $X,Y$ is a confusable pair, then so is $X\uplus Z,\,Y\uplus Z$ for every
multiset $Z$, and so is $X\setminus C,\,Y\setminus C$ with $C=X\cap Y$.
Consequently $D$ separates multisets of size $n$ if and only if $n<f(D)$, and a
smallest confusable pair has no point in common.
\end{lemma}

\begin{proof}
Projection respects sums, $\pi_d(X\uplus Z)=\pi_dX\uplus\pi_dZ$, and multisets
of numbers can be cancelled.
\end{proof}

Adding directions can only destroy confusable pairs, so $f(D')\le f(D)$
whenever $D'$ is part of $D$.  With Lemma~\ref{lem:pad} this turns bounds on
$f$ into bounds on $m_k$: for all $m$ and $n$,
\begin{equation}
\label{eq:flip}
  m_k(n)>m\iff\text{every list $D$ of $m$ directions has $f(D)\le n$,}
\end{equation}
since a list of fewer than $m$ directions is part of one of $m$ directions.
Lemma~\ref{lem:pad} also shows that a list that separates multisets of size
$n$ separates every smaller size, so $m_k$ is nondecreasing in $n$.

\paragraph{Signed multisets.}  It is easier to work with the difference of a
confusable pair.  A \emph{signed multiset} is a function $\nu:\R^k\to\Z$ with
finite support; we write $\pt{x}$ for the function that is $1$ at $x$ and $0$
elsewhere, and $\norm{\nu}:=\sum_x|\nu(x)|$.  For a direction $d$, let
$(\pi_d\nu)(t):=\sum_{x\cdot d=t}\nu(x)$ for $t\in\R$.  We say that $\nu$ is
\emph{silent} for $D$ if $\pi_d\nu=0$ for every $d$ in $D$.  For multisets,
$\pi_d(X-Y)=0$ says exactly that $X$ and $Y$ have the same projection on $d$.

\begin{lemma}
\label{lem:silent}
If $D$ contains at least one direction, then
\[
  f(D)=\min\{\tfrac12\norm{\nu} : \nu\neq0 \text{ is silent for } D\}.
\]
\end{lemma}

\begin{proof}
A confusable pair $X,Y$ of size $s$ gives the silent $\nu=X-Y\neq0$ with
$\norm{\nu}\le2s$.  Conversely, if $\nu\neq0$ is silent, its total charge is
$\sum_t(\pi_d\nu)(t)=0$ for a direction $d$ of $D$, so its positive and negative parts are multisets of
size $\frac12\norm\nu$, and they form a confusable pair.
\end{proof}

The \emph{convolution} of signed multisets is
$(\nu*\eta)(z):=\sum_{x+y=z}\nu(x)\eta(y)$, so that
$\pt{x}*\pt{y}=\pt{x+y}$.

\begin{lemma}
\label{lem:conv}
For signed multisets $\nu,\eta$ and a direction $d$:
\begin{enumerate}[label=(\roman*)]
\item $\pi_d(\nu*\eta)=(\pi_d\nu)*(\pi_d\eta)$, so $\nu*\eta$ is silent for
  every direction for which $\nu$ or $\eta$ is silent;
\item $\norm{\nu*\eta}\le\norm{\nu}\,\norm{\eta}$;
\item if $\nu\neq0$ and $\eta\neq0$ then $\nu*\eta\neq0$.
\end{enumerate}
\end{lemma}

\begin{proof}
(i) and (ii) are immediate from the definition.  For (iii), choose $\ell$ such
that $x\mapsto\ell\cdot x$ is injective on the supports of $\nu$ and $\eta$,
and let $x^*$, $y^*$ maximise $\ell\cdot x$ on them.  Then $x^*+y^*$ arises
only as $x^*+y^*$, so $(\nu*\eta)(x^*+y^*)=\nu(x^*)\eta(y^*)\neq0$.
\end{proof}

\begin{lemma}[Product lemma {\cite[Lemma~4.1]{MPS}}]
\label{lem:product}
For any lists $D_1$, $D_2$, the concatenated list $D_1\cup D_2$ has
$f(D_1\cup D_2)\le2f(D_1)f(D_2)$.
\end{lemma}

\begin{proof}
Take silent $\nu$, $\eta$ for $D_1$, $D_2$ with $\norm\nu=2f(D_1)$ and
$\norm\eta=2f(D_2)$ (Lemma~\ref{lem:silent}).  By Lemma~\ref{lem:conv},
$\nu*\eta$ is nonzero, silent for $D_1\cup D_2$, and has norm at most
$4f(D_1)f(D_2)$.
\end{proof}

\cite{MPS} state this for $k=2$; the proof works in every dimension.  For a
vector $v\neq0$, the \emph{dipole} $\pt0-\pt v$ is silent for exactly the
directions orthogonal to $v$.  Every $k-1$ directions are orthogonal to some
$v\neq0$, so a list of at most $k-1$ directions has $f=1$, and
Lemma~\ref{lem:product} gives the following.

\begin{lemma}[The hypercube {\cite[Theorem~5]{DWTLB}}]
\label{lem:hypercube}
Every list $D$ of $m$ directions in $\R^k$ has
$f(D)\le2^{\lceil m/(k-1)\rceil-1}$.
\end{lemma}

The corresponding silent signed multiset is
$(\pt0-\pt{v_1})*\cdots*(\pt0-\pt{v_M})=\sum_I(-1)^{|I|}\pt{v(I)}$, with
$M=\lceil m/(k-1)\rceil$, $v_j$ orthogonal to the $j$-th group of $k-1$
directions, and $v(I)=\sum_{j\in I}v_j$: the vertices of a parallelepiped,
coloured by parity.  The lengths of the $v_j$ are free, which
Section~\ref{sec:lower} exploits.

\begin{lemma}[One point at a time {\cite[Proposition~48 and Corollary~49]{HSW}}]
\label{lem:linear}
If $D$ is in general position (Definition~\ref{def:gp}), then
$f(D)\ge\lceil m/(k-1)\rceil$.  So
$(k-1)n+1$ directions in general position separate multisets of size $n$.
\end{lemma}

\begin{proof}
Let $X,Y$ be a smallest confusable pair, and fix $a\in X$.  For each direction
$d_j$ some $b\in Y$ has $b\cdot d_j=a\cdot d_j$.  Since $a\notin Y$
(Lemma~\ref{lem:pad}), each $b\in Y$ serves only the directions orthogonal to
$b-a\neq0$, of which there are at most $k-1$ (Lemma~\ref{lem:gp}(c)).  So
$|Y|\ge m/(k-1)$.
\end{proof}

Write two multisets of size $n$
as $X=\ms{x_1,\dots,x_n}$ and $Y=\ms{y_1,\dots,y_n}$.  They have the same
projection on $d_j$ exactly when, for some permutation $\rho_j$ of
$\{1,\dots,n\}$, $x_a\cdot d_j=y_{\rho_j(a)}\cdot d_j$ for every $a$.  For each
list $\rho=(\rho_1,\dots,\rho_m)$ of permutations these are linear equations in
the $2n$ unknown points; let $M_\rho(D)$ be their matrix, and
$K_\rho(D)\subseteq(\R^k)^{2n}$ the space of their solutions.  For a
permutation $\sigma$, let $Z_\sigma$ be the subspace on which
$x_a=y_{\sigma(a)}$ for every $a$, so that $X=Y$.

\begin{lemma}[A criterion for separation {\cite[Proposition~3.8]{BHS}}]
\label{lem:criterion}
A list $D$ separates multisets of size $n$ if and only if, for every list
$\rho$ of permutations, $K_\rho(D)$ lies in $Z_\sigma$ for some
permutation~$\sigma$.  Moreover, $K_\rho(D)$ lies in $Z_\sigma$ exactly when
adding the equations of $Z_\sigma$ to $M_\rho(D)$ does not raise its rank, so
whether $D$ separates multisets of size $n$ can be decided by finitely many
rank computations.
\end{lemma}

\begin{proof}
The multisets $X$ and $Y$ have the same projection on every direction exactly
when $(x_1,\dots,x_n,y_1,\dots,y_n)$ lies in some $K_\rho(D)$, and $X=Y$
exactly when it lies in some $Z_\sigma$.  So $D$ separates multisets of size
$n$ if and only if every $K_\rho(D)$ lies in the union of the subspaces
$Z_\sigma$, and a vector space that lies in a finite union of subspaces lies in
one of them.  Adding the equations of $Z_\sigma$ to $M_\rho(D)$ cuts its kernel
down to $K_\rho(D)\cap Z_\sigma$, which is all of $K_\rho(D)$ exactly when the
rank does not go up.  There are finitely many $\rho$ and $\sigma$.
\end{proof}

The text of \cite{BHS} states the lemma for lists that begin with the
coordinate axes.  For rational directions the rank computations can be done in
exact arithmetic.  Relabelling $Y$ lets us take $\rho_1$ to be the
identity, which leaves $(n!)^{m-1}$ lists $\rho$, so the test is practical only
for small $n$ and $m$ (Remark~\ref{rem:small}).  This test was used
numerically in \cite[Section~4]{DWTLB}.

\paragraph{Witnesses.}  Several proofs below show that a property holds for
generic lists of directions (Definition~\ref{def:ae}) by finding a polynomial
that vanishes at every list that fails, and a \emph{witness}: one list, however
special, at which the polynomial is not zero.  Dym and Gortler turned this idea
into a count, the \emph{finite witness theorem} \cite[Theorem~2.7]{DG}, which
Amir et al.\ \cite{AGARD} use for multisets and Dym et al.\
\cite[Theorem~7]{DWTLB} for sorted projections; we use it in
Remark~\ref{rem:compress}.

\section{Lower bounds}
\label{sec:lower}

\begin{proof}[Proof of Theorem~\ref{thm:three}]
Split $D$ into two groups $G_1$, $G_2$ of $k-1$ directions and one more
direction $d$, and choose nonzero $v_1\perp G_1$ and $v_2\perp G_2$.

If $v_1$ and $v_2$ are parallel, the dipole $\pt0-\pt{v_1}$ is silent for
$G_1$ and $G_2$; with $u\neq0$ orthogonal to $d$, the signed multiset
$(\pt0-\pt{v_1})*(\pt0-\pt u)$ is nonzero, silent for $D$ and of norm at most
$4$, so $f(D)\le2$.

Otherwise the plane spanned by $v_1$, $v_2$ meets $d^\perp$ in some
$w=\alpha v_1+\beta v_2\neq0$.  If $\beta=0$ then $v_1\perp d$, and
$(\pt0-\pt{v_1})*(\pt0-\pt{v_2})$ gives $f(D)\le2$; likewise if $\alpha=0$.
If $\alpha\beta\neq0$, put $w_1=\alpha v_1$ and $w_2=\beta v_2$, so that
$w_1+w_2\perp d$, and let
\[
  \nu=(\pt0-\pt{w_1})*(\pt0-\pt{w_2})*(\pt0-\pt{w_1+w_2}).
\]
It is nonzero and silent for $D$ (Lemma~\ref{lem:conv}).  Of the eight terms in
its expansion, $+\pt{w_1+w_2}$ and $-\pt{w_1+w_2}$ cancel, so
$\norm\nu\le6$ and $f(D)\le3$.

(a) By the above, and since fewer directions can only do worse,
$m_k(3)\ge2k$; Lemma~\ref{lem:linear} gives $m_k(3)\le3(k-1)+1$.
(b) Let $D$ be any list of $q(k-1)+1=(2k-1)+(q-2)(k-1)$ directions.  It splits
into $2k-1$ directions, with $f\le3$ by the above, and $q-2$ groups of $k-1$,
with $f=1$, so Lemma~\ref{lem:product} and the hypothesis $3\cdot2^{q-2}\le n$
give $f(D)\le3\cdot2^{q-2}\le n$.  By \eqref{eq:flip}, $m_k(n)\ge q(k-1)+2$.
The other hypothesis, $n<2^q$, is needed only to rewrite this: it gives
$\lfloor\log_2n\rfloor+1\le q$, so
$q(k-1)+2\ge(k-1)\bigl(\lfloor\log_2n\rfloor+1\bigr)+2$.
\end{proof}

For $k=2$ the third factor of $\nu$ is the triangle of \cite[Section~4]{MPS}.
The proof uses the freedom in the lengths of the $v_j$ noted after
Lemma~\ref{lem:hypercube}: with $2k-1$ directions the third dipole can be
moved into the plane of the first two.

\begin{remark}[Small cases by computer]
\label{rem:small}
By Lemma~\ref{lem:criterion}, whether a given list separates multisets of size
$n$ can be decided by a finite computation, in exact arithmetic when the
directions are rational.  It shows that
$e_1,e_2,e_3,(1,2,4),(3,-1,2),(-2,5,1)$ separates multisets of size $3$ in
$\R^3$, and that
$e_1,\dots,e_4,(5,2,-8,3),(-1,7,4,9),(8,-3,6,-5),(4,9,-2,7)$ does so in $\R^4$.
So $m_3(3)=6$ and $m_4(3)=8$.  In the plane, $(1,0),(0,1),(1,1),(1,3)$
separates multisets of size $4$, and adding $(3,-2)$ gives size $5$; with
Theorem~\ref{thm:plane} below this gives $m_2(4)=4$ and $m_2(5)=5$.  For sets
rather than multisets, \cite[Propositions~1 and~2]{BL90} already show that any
$m$ directions in the plane that are not \emph{affinely rational} (not mapped by
a linear map to directions at rational multiples of $\pi$) determine every set
of at most $m$ points.  These checks use \texttt{locate\_exact.py} in the accompanying code, which runs the
test in integer arithmetic.

Dym et al.\ \cite[Section~4]{DWTLB} found lists that separate in the cases
$(n,k,m)=(4,2,4)$, $(5,2,5)$ and $(3,4,8)$ with a floating-point version of the
test, and the lists they print pass the exact test too.  They also report a list
for $(3,3,6)$, but the list as printed does not separate: three of its
directions, $e_1$, $e_3$ and $(0.21,0,0.45)$, lie in one plane, so the
construction of Theorem~\ref{thm:three} fools it.  Presumably an entry was
rounded for display.  They also observed that random lists of $2k-1$ directions
fail for $k=2,3,4$, which Theorem~\ref{thm:three} proves for every list and
every $k$.

For $k\ge3$, general position is not enough at $m=2k$.  The list
$e_1$, $e_2$, $e_3$, $(1,1,1)$, $(1,2,3)$, $(2,-1,1)$ is in general position but does not
separate multisets of size $3$: it splits into the pairs $\{e_1,e_3\}$,
$\{e_2,(1,1,1)\}$ and $\{(1,2,3),(2,-1,1)\}$, whose common normals $(0,1,0)$,
$(1,0,-1)$ and $(1,1,-1)$ are linearly dependent, and that is all the
construction above needs.  Such lists exist in every dimension $k\ge3$: take
linearly independent $v_1,v_2$ and choose $k-1$ directions orthogonal to
$v_1$, $k-1$ orthogonal to $v_2$ and $2$ orthogonal to $v_1+v_2$.  The
construction above fools every such list.  For almost every such choice
(Definition~\ref{def:ae}) the list is in general position, because any $k$ of
the directions include two from different groups,
and choosing those two, each in its own hyperplane, outside
$v_1^\perp\cap v_2^\perp$, and the others as a basis of $v_1^\perp\cap
v_2^\perp$, gives a basis of $\R^k$; so no $k\times k$ determinant vanishes
identically.
\end{remark}

\paragraph{Specialising to $k=2$:}
In the plane the construction of Theorem~\ref{thm:three} bounds $f(D)$ above by a quantity of order  $6^{m/3}$. The best construction of
\cite{MPS}, with nine directions, gives $f(D)<O(1.7996^m)$.  
The better bound on $f(D)$ claimed by Theorem~\ref{thm:plane} relies on results concerning groups of four and five directions in the plane that the next two Lemmas provide.

\begin{lemma}[Bianchi and Longinetti {\cite[p.~234 and Fig.~2]{BL90}}]
\label{lem:four}
Every list $D$ of $4$ directions in the plane has $f(D)\le5$.
\end{lemma}

\begin{proof}
If two directions are parallel there are at most three different directions,
and $f(D)\le3$ by Theorem~\ref{thm:three}.  Otherwise, since the projection of
$X$ on $Ad$ is the projection of $A^{\mathsf T}X$ on $d$, an invertible linear
map lets us assume that the directions are $(1,0)$, $(0,1)$, $(1,1)$ and
$(1,\lambda)$ with $\lambda\neq0,1$.  Let
\begin{align*}
  X&=\ms{(0,0),\,(-\lambda(\lambda-1),(\lambda-1)^2),\,(-\lambda^2,1),\,
       (-\lambda,\lambda^2-\lambda+1),\,(-\lambda(2\lambda-1),\lambda^2)},\\
  Y&=\ms{(0,(\lambda-1)^2),\,(-\lambda(\lambda-1),0),\,(-\lambda^2,\lambda^2),\,
       (-\lambda,1),\,(-\lambda(2\lambda-1),\lambda^2-\lambda+1)}.
\end{align*}
With $\mu=1-\lambda$, the projections of the points, in the order listed, are
\begin{center}
\renewcommand{\arraystretch}{1.15}
\begin{tabular}{lll}
  direction & $X$ & $Y$\\\hline
  $(1,0)$ & $0,\ \lambda\mu,\ -\lambda^2,\ -\lambda,\ \lambda(1-2\lambda)$
          & $0,\ \lambda\mu,\ -\lambda^2,\ -\lambda,\ \lambda(1-2\lambda)$\\
  $(0,1)$ & $0,\ \mu^2,\ 1,\ \lambda^2-\lambda+1,\ \lambda^2$
          & $\mu^2,\ 0,\ \lambda^2,\ 1,\ \lambda^2-\lambda+1$\\
  $(1,1)$ & $0,\ \mu,\ 1-\lambda^2,\ \mu^2,\ \lambda\mu$
          & $\mu^2,\ \lambda\mu,\ 0,\ \mu,\ 1-\lambda^2$\\
  $(1,\lambda)$ & $0,\ \lambda\mu(1+\mu),\ \lambda\mu,\ -\lambda^2\mu,\ \lambda\mu^2$
          & $\lambda\mu^2,\ \lambda\mu,\ -\lambda^2\mu,\ 0,\ \lambda\mu(1+\mu)$
\end{tabular}
\end{center}
and in each row the two lists are rearrangements of each other.  The multisets
differ, because $(0,0)\in X$ while the points of $Y$ with first coordinate $0$
are $(0,\mu^2)$ and, if $\lambda=\frac12$, $(0,\frac34)$.
\end{proof}

\cite{BL90} give this example for sets, to show that their Proposition~3 needs
its hypothesis $m>4h^2+11h$ (here $h=1$ and $m=4$).  That proposition says that
$m$ directions that are not affinely rational determine every set of $m+h$
points; it is the only result of \cite{BL90} that \cite{MPS} cite, as a lower
bound $k+\Omega(\sqrt k)$ for sets.

\begin{lemma}
\label{lem:five}
Every list $D$ of $5$ directions in the plane has $f(D)\le6$.
\end{lemma}

\begin{proof}
If two directions are parallel there are at most four different directions,
and $f(D)\le5$ by Lemma~\ref{lem:four}.  Otherwise an invertible linear map lets
us assume that the directions are $(1,0)$, $(0,1)$, $(1,1)$, $(1,a)$ and
$(1,b)$, with $a,b\notin\{0,1\}$ and $a\neq b$.  Let $X=\ms{x_1,\dots,x_6}$ with
\begin{gather*}
  x_1=(ab-a+b,\ a-b-1),\quad x_2=(b,\,-1),\quad x_3=(ab,\,-b),\\
  x_4=(b-a,\ a-b),\quad x_5=(ab-a,\ 0),\quad x_6=(0,\ a-1),
\end{gather*}
let $c=x_1$, and let $Y=\ms{c-x_1,\dots,c-x_6}$ be the reflection of $X$ in the
point $c/2$.  For each direction $d$ the projections $x_i\cdot d$ fall into
three pairs, each with sum $c\cdot d$:
\begin{center}
\renewcommand{\arraystretch}{1.15}
\begin{tabular}{@{}l@{\quad}l@{\quad}l@{\quad}l@{}}
  $d$ & $x_1\cdot d,\dots,x_6\cdot d$ & $c\cdot d$ & pairs\\\hline
  $(1,0)$ & $ab-a+b,\ b,\ ab,\ b-a,\ ab-a,\ 0$ & $ab-a+b$ & $16,\ 25,\ 34$\\
  $(0,1)$ & $a-b-1,\ -1,\ -b,\ a-b,\ 0,\ a-1$ & $a-b-1$ & $15,\ 24,\ 36$\\
  $(1,1)$ & $ab-1,\ b-1,\ ab-b,\ 0,\ ab-a,\ a-1$ & $ab-1$ & $14,\ 23,\ 56$\\
  $(1,a)$ & $a^2-2a+b,\ b-a,\ 0,\ a^2-ab-a+b,\ ab-a,\ a^2-a$ & $a^2-2a+b$
          & $13,\ 26,\ 45$\\
  $(1,b)$ & $2ab-a-b^2,\ 0,\ ab-b^2,\ ab-a+b-b^2,\ ab-a,\ ab-b$ & $2ab-a-b^2$
          & $12,\ 35,\ 46$
\end{tabular}
\end{center}
(``$16$'' means $x_1\cdot d+x_6\cdot d=c\cdot d$, and so on).  So the projection
of $X$ on $d$ is unchanged by $t\mapsto c\cdot d-t$, which takes it to the
projection of $Y$.  The multisets differ.  Every point of $X$ or $Y$ has first
coordinate $0$, $ab-a+b$, $b$, $ab$, $b-a$ or $ab-a$, and the last four are not
$0$.  The points of $X$ with first coordinate $0$ are $x_6=(0,a-1)$ and, if
$ab-a+b=0$, also $x_1$; those of $Y$ are $c-x_1=(0,0)$ and, if $ab-a+b=0$,
also $c-x_6=(0,-b)$.  Their second coordinates add up to $a-1$ and $0$, or to
$2a-b-2$ and $-b$, and these differ because $a\neq1$.
\end{proof}

We found Lemma~\ref{lem:five} by a computer search over the matchings that a
pair with equal projections induces between its two multisets
(Remark~\ref{rem:six-seven}); \texttt{pair6.py} in the accompanying code checks
the lemma.  The
bound $6$
is attained: the list $(1,0)$, $(0,1)$, $(1,1)$, $(1,3)$, $(3,-2)$ separates multisets of
size $5$ (Remark~\ref{rem:small}).  For single chosen lists of five lattice
directions the same pair was known \cite[after Theorem~2.2]{AL15},
\cite[Table~3.1]{Ghi19}, and \cite[Theorem~3.8.4]{Ghi19} combines it with
specially constructed directions to give pairs of size $O(12^{m/5})$ for those
directions.  We know of no earlier statement for every list of five
directions.

\begin{proof}[Proof of Theorem~\ref{thm:plane}]
[This proof uses the values 
$(c_0,\dots,c_4)=(1,2,4,6,10)$ 
and the decomposition $m=5 t + m$
from the statement of the theorem.]
Split $D$ into $t$ lists $D_1,\ldots,D_t$ of five directions each, in addition to a single remainder group $R$ of $r$ directions if $r>0$.
Each $D_i$ has $f(D_i)\le6$
(Lemma~\ref{lem:five}).  A remainder group of four directions has $f(R)\le5\equiv c_4/2$ (Lemma~\ref{lem:four}), one of
three has $f(R)\le3\equiv c_3/2$ (Theorem~\ref{thm:three}), one of two has $f(R)\le2\equiv c_2/2$
(Lemma~\ref{lem:hypercube}), and one of one has $f(R)=1\equiv c_1/2$.   When $r=0$ the product lemma
(Lemma~\ref{lem:product}), acting just on the $t$ groups of five, says that $f(D)\le 2^{t-1} 6^t$.  When $r>0$ the product lemma picks up an extra factor of 2 and an extra $f(R)$  from the remainder group, so  
$f(D)\le 2^t 6^t f(R)=2^t 6^t (c_r/2)$.  These two bounds can be combined in the single expression
$f(D)\le \frac12c_r12^t$.  We do not believe any other grouping strategy among 
groups of one to five directions does better.  As $c_r\le12$,
$f(D)\le\frac12\,12^{\lceil m/5\rceil}$, so if $D$ separates multisets of size
$n$ then $\lceil m/5\rceil>\log_{12}(2n)$ and $m>5\log_{12}(2n)-5$.
\end{proof}

\begin{remark}[Six and seven directions in the plane]
\label{rem:six-seven}
An exhaustive search over matching graphs, described in the rest of this
remark, shows that
generic lists (Definition~\ref{def:ae}) of $4$, $5$, $6$ and $7$ directions in
the plane separate multisets of size $4$, $5$, $11$ and $14$ respectively.
Since every list of five directions has $f(D)\le6$ (Lemma~\ref{lem:five}) and
every list of six has $f(D)\le12$ (Theorem~\ref{thm:plane}), this gives
$m_2(n)=6$ for $6\le n\le11$ and $m_2(n)=7$ for
$12\le n\le14$.\footnote{After reading the first public version of this paper,
\cite{ZM26} told us that they had shown $m_2(6)=6$, which this
remark confirms.  Apart from their suggestion acknowledged after
Theorem~\ref{thm:plane}, which led to Lemma~\ref{lem:five}, we have used
nothing from their message.  The search follows up our remark in that version
that six generic directions are not fooled by a $6+6$ pair \cite{LS26}, and it
goes further: almost every list of six directions separates multisets of size
$11$, and $m_2(n)$ is settled for every $n\le14$.}  So in the plane, for every
$n\le14$ (for $n\le3$ by the linear bound), almost every list
of $m_2(n)$ directions separates multisets of size $n$.  Almost every list of six directions has $f(D)=12$, so
groups of six directions give at best the base $24^{1/6}\approx1.698$ in
Theorem~\ref{thm:plane}, worse than $12^{1/5}$.

Let $X=\ms{x_1,\dots,x_s}$ and $Y=\ms{y_1,\dots,y_s}$ be a smallest confusable
pair for a list of $m$ directions, and let $d_j$ match $x_i$ with
$y_{\sigma_j(i)}$ (Definition~\ref{def:matching}).  Two non-parallel
directions never match the same two points, which would then coincide, against
Lemma~\ref{lem:pad}; so the $\sigma_j$ form an $m\times s$ Latin rectangle.
For fixed $\sigma_j$ the conditions $(x_i-y_{\sigma_j(i)})\cdot d_j=0$ are
linear, and their solutions form a space $S$ that contains the translations.
For every Latin rectangle, up to relabelling, exchanging $X$ and $Y$ and
reordering the directions, the search decides whether it \emph{works}: whether,
for generic lists, $S$ contains a pair with $X\neq Y$ and no point in common.
Ranks are computed modulo a large prime at random directions, which can only
overestimate $\dim S$, so no rectangle is lost, and every rectangle that
survives is decided exactly, in rational functions of the coordinates of the
directions.  Three facts prune the search.
\begin{enumerate}[label=(\roman*),leftmargin=*,itemsep=1pt]
\item If a rectangle works, so does every rectangle formed by some of its rows,
  and deleting a direction raises the generic $\dim S$ by at least one.
  Otherwise the equations of $d_j$ would hold on the solutions for the other
  directions, which do not depend on $d_j$, for generic and hence for all
  $d_j$, and then $x_i=y_{\sigma_j(i)}$ for every $i$.  So any $r$ of its $m$
  rows have generic $\dim S\ge3+m-r$.
\item For disjoint sets $A$ and $B$ of directions, generic
  $\dim S(A\cup B)\le c(A)+c(B)$, where $c(A)$ is the number of connected
  components of the matching graph in the colours of $A$.  Indeed, giving every
  direction of $A$ one value $d$ and every direction of $B$ another value $e$
  can only lower the rank, and it leaves just the conditions that $x\cdot d$ be
  constant on each component for $A$, and $x\cdot e$ on each component for
  $B$.  For seven directions, (i) and~(ii) show that some four directions are
  such that the matching graph in the colours of any three of them is
  disconnected, and in the colours of any two has at least three components.
  Take three directions $T$ with $c(T)$ smallest.  If $c(T)\le2$, the other
  four will do, by (i) and~(ii) for $T$ with three or two of them.  Otherwise
  any four will, since adding colours only merges components.  The search
  starts from such four.
\item If a component of the matching graph in some of the colours has, by its
  own equations, only the translations as solutions for generic lists, then
  its points coincide, and $X$ and $Y$ have a point in common.  A rank modulo
  the prime that leaves only the translations proves this.
\end{enumerate}
The search is \texttt{latin\_search/} in the accompanying code, where
\texttt{run\_claim.sh} with the argument \texttt{m4}, \texttt{m5}, \texttt{m6}
or \texttt{m7} reproduces each of the four claims and its \texttt{README.md} explains the steps.  The search for
seven directions and $14$ points took about $35$ CPU-hours.
\end{remark}

\section{Upper bounds from geometric rigidity}
\label{sec:rigidity}

Throughout this section and the next, $X,Y$ is a confusable pair for a list
$D=(d_1,\dots,d_m)$, of some size $s$, with no point in common (for instance a
smallest one, Lemma~\ref{lem:pad}).  We treat its $2s$ points as $2s$ separate
\emph{vertices}, even when two points of $X$ coincide.

\begin{definition}[Matching graph]
\label{def:matching}
For each $j$ the multisets $\pi_{d_j}X$ and $\pi_{d_j}Y$ are equal, so there is
a perfect matching $E_j$ between the vertices of $X$ and those of $Y$ with
$(a-b)\cdot d_j=0$ for every pair $\{a,b\}\in E_j$; we fix one.  The
\emph{matching graph} of the pair has these $2s$ vertices and, for each $j$, the
edges of $E_j$ in \emph{colour}~$j$.
\end{definition}

So every colour class is a perfect matching, and every vertex has exactly one
edge of each colour.  An edge of colour $j$ forces $(a-b)\cdot d_j=0$, but not
conversely: when several points have the same projection on $d_j$, the matching
we fix may pair them differently.

\begin{definition}[Pinning pattern {\cite[Definition~8]{HS26}}]
\label{def:pattern}
A \emph{pinning pattern} in an edge-coloured multigraph is a set $W$ of at least
two vertices together with $k$ spanning trees $T_1,\dots,T_k$ of $W$, made of
edges with both ends in $W$, such that no colour is used by two different trees.
\end{definition}

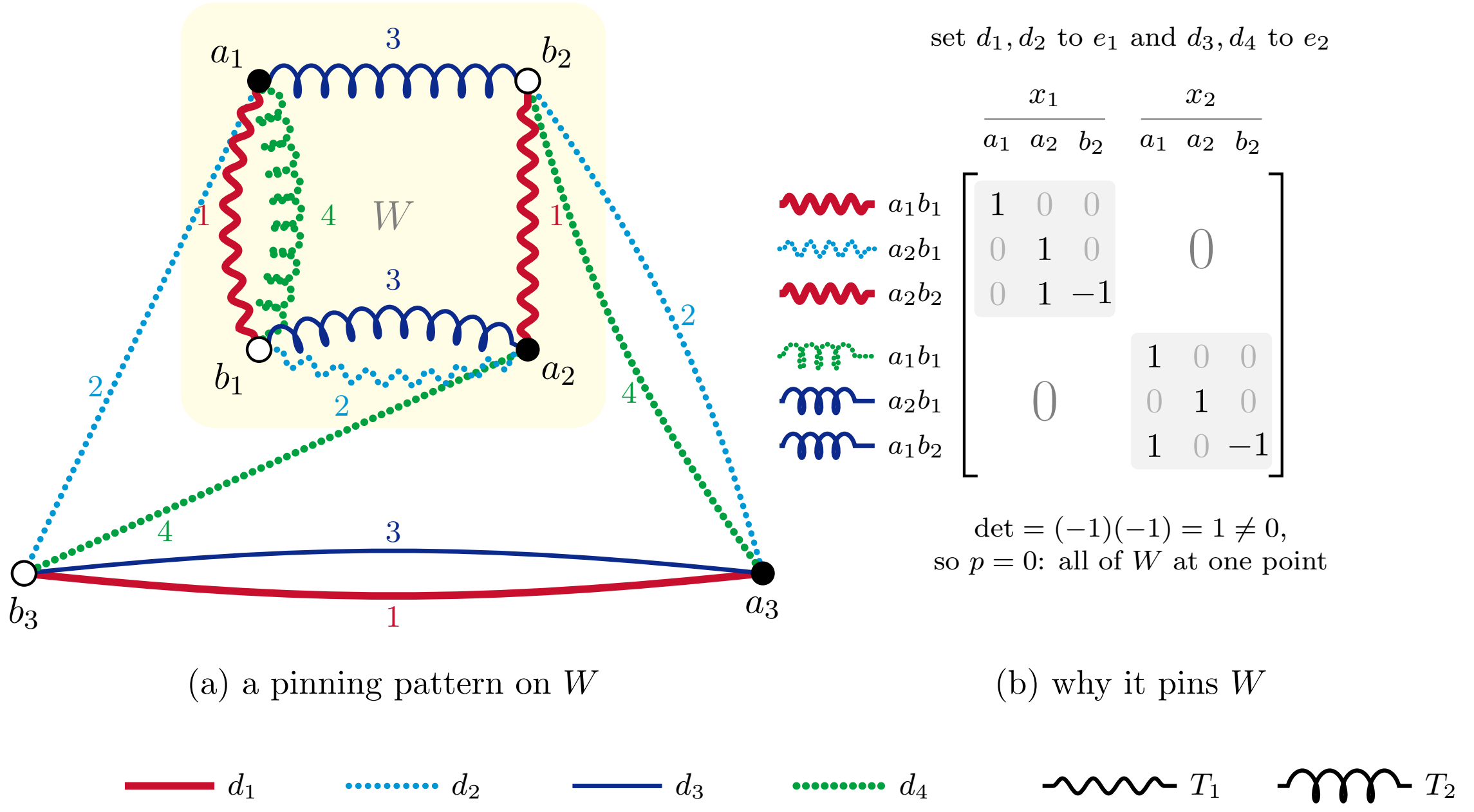
\begin{figure}[htbp]
\centering
\begin{tikzpicture}[scale=1.15,
    xm/.style={circle,fill=black,inner sep=0pt,minimum size=7pt},
    ym/.style={circle,draw=black,thick,fill=white,inner sep=0pt,minimum size=7pt},
    lab/.style={font=\footnotesize,inner sep=2pt}]
  \begin{scope}
    \fill[yellow!12,rounded corners=10pt] (-0.7,-0.7) rectangle (3.1,3.1);
    \node[gray,font=\large] at (1.2,1.2) {$W$};
    \coordinate (a1) at (0,2.4);  \coordinate (b2) at (2.4,2.4);
    \coordinate (b1) at (0,0);    \coordinate (a2) at (2.4,0);
    \coordinate (b3) at (-2.1,-2.0); \coordinate (a3) at (4.5,-2.0);
    \draw[bat1] (b3) to[bend right=6] node[lab,below,text=batA] {$1$} (a3);
    \draw[bat3] (b3) to[bend left=6]  node[lab,above,text=batC] {$3$} (a3);
    \draw[bat2] (a1) -- node[lab,left,text=batB,pos=0.62] {$2$} (b3);
    \draw[bat4] (a2) -- node[lab,below,text=batD,pos=0.72] {$4$} (b3);
    \draw[bat2] (b2) to[bend left=10]  node[lab,right,text=batB] {$2$} (a3);
    \draw[bat4] (b2) to[bend right=10] node[lab,left,text=batD,pos=0.62] {$4$} (a3);
    \draw[bat1,photon] (a1) to[bend right=20] node[lab,left=3pt,text=batA] {$1$} (b1);
    \draw[bat1,photon] (a2) -- node[lab,right=3pt,text=batA] {$1$} (b2);
    \draw[bat2,photon] (b1) to[bend right=20] node[lab,below=4pt,text=batB,pos=0.3] {$2$} (a2);
    \draw[bat4,spring] (a1) to[bend left=20] node[lab,right=7pt,text=batD] {$4$} (b1);
    \draw[bat3,spring] (a1) -- node[lab,above=7pt,text=batC] {$3$} (b2);
    \draw[bat3,spring] (b1) to[bend left=20] node[lab,above=7pt,text=batC] {$3$} (a2);
    \node[xm] at (a1) {}; \node[xm] at (a2) {}; \node[xm] at (a3) {};
    \node[ym] at (b1) {}; \node[ym] at (b2) {}; \node[ym] at (b3) {};
    \node[above left]  at (a1) {$a_1$};  \node[above right] at (b2) {$b_2$};
    \node[below left]  at (b1) {$b_1$};  \node[below right] at (a2) {$a_2$};
    \node[below=3pt]   at (b3) {$b_3$};  \node[below=3pt]   at (a3) {$a_3$};
    \node[font=\small] at (1.2,-3.0) {(a) a pinning pattern on $W$};
  \end{scope}
  \begin{scope}[xshift=6.6cm,yshift=1.3cm,
      rowplain/.style={},
      rowdense/.style={line width=1.4pt,dash pattern=on 0pt off 1.6pt}]
    \fill[gray!10,rounded corners=2pt] (-0.21,0.21) rectangle (1.05,-1.01);
    \fill[gray!10,rounded corners=2pt] (1.19,-1.15) rectangle (2.45,-2.37);
    \draw[thick] (-0.18,0.27) -- (-0.30,0.27) -- (-0.30,-2.43) -- (-0.18,-2.43);
    \draw[thick] (2.42,0.27) -- (2.54,0.27) -- (2.54,-2.43) -- (2.42,-2.43);
    \foreach \x/\v in {0/a_1,0.42/a_2,0.84/b_2,1.40/a_1,1.82/a_2,2.24/b_2}
      \node[font=\scriptsize] at (\x,0.55) {$\v$};
    \node[font=\footnotesize] at (0.42,0.95) {$x_1$};
    \node[font=\footnotesize] at (1.82,0.95) {$x_2$};
    \draw[gray] (-0.12,0.76) -- (0.96,0.76);
    \draw[gray] (1.28,0.76) -- (2.36,0.76);
    \node[anchor=south,align=center,font=\scriptsize] at (1.19,1.25)
      {set $d_1,d_2$ to $e_1$ and $d_3,d_4$ to $e_2$};
    \foreach \y/\s/\t/\v in {0/bat1/rowplain/a_1b_1,-0.40/bat2/rowdense/a_2b_1,
                             -0.80/bat1/rowplain/a_2b_2}{
      \draw[\s,\t,photon,decoration={segment length=6pt,pre length=1pt,
        post length=1pt}] (-1.95,\y) -- (-1.10,\y);
      \node[anchor=east,font=\scriptsize] at (-0.36,\y) {$\v$};
    }
    \foreach \y/\s/\t/\v in {-1.36/bat4/rowdense/a_1b_1,-1.76/bat3/rowplain/a_2b_1,
                             -2.16/bat3/rowplain/a_1b_2}{
      \draw[\s,\t,spring,decoration={amplitude=3.5pt,segment length=5pt,
        pre length=1pt,post length=1pt}] (-1.95,\y) -- (-1.10,\y);
      \node[anchor=east,font=\scriptsize] at (-0.36,\y) {$\v$};
    }
    \foreach \y/\e in {0/{1,0,0},-0.40/{0,1,0},-0.80/{0,1,-1}}{
      \foreach \c [count=\i from 0] in \e {
        \ifnum\c=0 \node[gray!60,font=\small] at (0.42*\i,\y) {$0$};
        \else \node[font=\small] at (0.42*\i,\y) {$\c$}; \fi
      }
    }
    \foreach \y/\e in {-1.36/{1,0,0},-1.76/{0,1,0},-2.16/{1,0,-1}}{
      \foreach \c [count=\i from 0] in \e {
        \ifnum\c=0 \node[gray!60,font=\small] at (1.40+0.42*\i,\y) {$0$};
        \else \node[font=\small] at (1.40+0.42*\i,\y) {$\c$}; \fi
      }
    }
    \node[gray,font=\LARGE] at (1.82,-0.40) {$0$};
    \node[gray,font=\LARGE] at (0.42,-1.76) {$0$};
    \node[font=\scriptsize,align=center] at (1.19,-3.05)
      {$\det=(-1)(-1)=1\neq0$,\\ so $p=0$: all of $W$ at one point};
    \node[font=\small] at (1.19,-4.3) {(b) why it pins $W$};
  \end{scope}
  \begin{scope}[yshift=-3.9cm,xshift=-1.2cm,font=\footnotesize]
    \foreach \b/\x in {1/0,2/2.0,3/4.0,4/6.0} {
      \draw[bat\b] (\x,0) -- ++(0.8,0);
      \node[right] at (\x+0.8,0) {$d_{\b}$};
    }
    \draw[line width=1.2pt,photon] (8.2,0) -- ++(1.2,0);
    \node[right] at (9.4,0) {$T_1$};
    \draw[line width=1.2pt,spring] (10.3,0) -- ++(1.2,0);
    \node[right] at (11.5,0) {$T_2$};
  \end{scope}
\end{tikzpicture}
\caption[A pinning pattern]{A pinning pattern for $k=2$, drawn as an abstract
coloured graph (the positions of the vertices mean nothing).  (a)~A matching
graph on six points, $a_1,a_2,a_3$ from $X$ (filled) and $b_1,b_2,b_3$ from $Y$
(open), for four directions.  The edges of colour $j$, a perfect matching, have
the colour and line style of $d_j$ in the key.  So every point has exactly one
edge of each colour, and edges that meet at a point always have different
colours.  On the shaded set $W=\{a_1,a_2,b_1,b_2\}$, the tree $T_1$ is drawn with
wiggles \photonsample{} and uses colours $1$ and $2$; the tree $T_2$ is drawn
with springs \springsample{} and uses colours $3$ and $4$.  Each tree needs two
colours: colour~$1$ pairs $W$ into $\{a_1,b_1\}$ and $\{a_2,b_2\}$, and
colour~$2$ joins the two pairs, and likewise for colours $3$ and $4$; since edges
that meet have different colours, the two colours of each tree alternate along
it.  The trees share vertices, and even the pairs $a_1b_1$ and $a_2b_1$, but no
colour.  Straight edges are used by neither tree.  (b)~Why such a pattern pins
$W$, as in the proof of Lemma~\ref{lem:rigidity}.  With $b_1$ at the origin, the
six tree edges give a square system $M(D)\,p=0$ in the coordinates of
$a_1,a_2,b_2$; the row of an edge of colour $j$ carries the first coordinate of
$d_j$ in the $x_1$ columns and its second in the $x_2$ columns.  Shown is $M$ at
the useless but legal directions $d_1=d_2=e_1$, $d_3=d_4=e_2$ (a witness, in
the sense of Section~\ref{sec:prelim}), which are consistent because no colour
is in both trees.  The rows of $T_1$ lose their
$x_2$ entries and those of $T_2$ their $x_1$ entries, so $M$ is block-diagonal,
and each block, a spanning tree taken outward from $b_1$, is triangular with
determinant $-1$.  So $P(D)=\det M(D)$ is $1$ here, and is not identically zero.
(For these trees, $P(D)=\det(d_1,d_4)\det(d_2,d_3)\det(d_1,d_3)$.)}
\label{fig:pattern}
\end{figure}

A tree may use several colours, and must do so unless $|W|=2$, because a colour
class is a matching.  Hair and Sahai call this a \emph{tree certificate}; for
$k=2$ it is the pair of colour-disjoint spanning trees of
\cite[Section~2.2]{MPS}.  Figure~\ref{fig:pattern} shows an example.
Section~\ref{sec:graphs} also needs a \emph{pinning pattern with $d$ trees},
defined in the same way with $d$ spanning trees in place of $k$.

\begin{lemma}[Rigidity {\cite[Lemma~7]{HS26}}]
\label{lem:rigidity}
If the matching graph of a confusable pair with no point in common contains a
pinning pattern, then the coordinates of the directions satisfy $P(D)=0$, where
$P$ is a polynomial with integer coefficients, depending only on the coloured
graph and the pattern, that is not identically zero.
\end{lemma}

\begin{proof}
Let the pattern be $T_1,\dots,T_k$ on $W$, fix $w_0\in W$, and treat the
positions $p_w\in\R^k$ of the other vertices of $W$ as unknowns, with
$p_{w_0}=0$.  The tree edges give $k(|W|-1)$ linear equations
$(p_a-p_b)\cdot d_j=0$, one for each edge $\{a,b\}$ of colour $j$ in one of the
trees, in $k(|W|-1)$ unknowns.  Write them as $M(D)\,p=0$ and put
$P(D):=\det M(D)$.

For a witness, set every direction used by $T_r$ equal to $e_r$; this is
consistent because no direction is used by two trees.  Then the equations of $T_r$ say that the ends of
each of its edges have equal $r$-th coordinates.  As $T_r$ spans $W$ and
$p_{w_0}=0$, all the $r$-th coordinates vanish.  This holds for every $r$, so
$p=0$ is the only solution, and $P\neq0$ there, so $P$ is not identically zero.

For the true directions, the true points, translated so that $w_0$ is at the
origin, solve $M(D)\,p=0$.  This solution is not $0$, because an edge of $T_1$ at
$w_0$ joins two different points.  So $P(D)=0$.
\end{proof}

For $k=2$ this is \cite[Lemmas~2.1--2.3]{MPS}.  Hair and Sahai prove it as
stated, for linear forms over any field, by the same evaluation.

\begin{corollary}
\label{cor:generic-no-pattern}
Fix $n$ and $m$.  For generic lists $D$ of $m$ directions, and so for almost
every list (Definition~\ref{def:ae}), no confusable pair of size at most $n$ with no point in common has a matching graph that
contains a pinning pattern.
\end{corollary}

\begin{proof}
Such a matching graph is one of finitely many coloured graphs on at most $2n$
vertices with $m$ colours, and a pinning pattern in it is one of finitely many
choices.  Each choice gives a polynomial that is not identically zero
(Lemma~\ref{lem:rigidity}), and $D$ lies in its zero set.  So every list that
fails is a zero of the product of these finitely many polynomials, which has
integer coefficients and is not identically zero; its zero set has measure zero
(Lemma~\ref{lem:generic}).
\end{proof}

\section{Upper bounds from graph theory}
\label{sec:graphs}

Theorem~\ref{thm:main} needs, besides Corollary~\ref{cor:generic-no-pattern},
a supply of pinning patterns: a matching graph with enough colours always
contains one.  Hair and Sahai prove such a lemma for list recovery of random
linear codes.  We state it in the form that their proof gives, and then tighten
it.

\begin{lemma}[{\cite[Lemma~6]{HS26}}]
\label{lem:hs}
Let $b\ge2$ be the number of vertices, $d\ge1$ be the number of spanning trees
required, and let $\beta>0$ be the matching density factor.  The number of
colours per tree is given by
\[
  T(b)\;:=\;\Bigl\lceil\frac{4\ln b}{\beta}\bigl(\ln\ln(4b)+\ln(2d)+1\bigr)
  \Bigr\rceil .
\]
Let $G$ be a loopless edge-coloured multigraph on $b$ vertices in which every
colour class is a matching with at least $\beta b$ edges.  If $G$ has at least
$d\,T(b)$ colours, then it contains a pinning pattern with $d$ trees.
\end{lemma}

Hair and Sahai state the lemma with a constant that they do not evaluate; the
form above is what their proof uses, padded so that its order of growth is
plain.  Their proof takes for $W$ a set of vertices
that maximises $e(W)/(|W|\ln|W|)$, where $e(W)$ counts the edges with both ends
in $W$.  By maximality, every split of $W$ into $A$ and $W\setminus A$ with
$|A|\le|W|/2$ is crossed by at least $(\beta m/\ln b)\ln(|W|/|A|)$ different
colours, since a matching puts at most $|A|$ edges across it.  They then choose
$d$ disjoint random sets of $T(b)$ colours each.  By their Lemma~5, a version of
the random-colour argument of \cite[Section~3.2]{MPS}, each set connects $W$
with probability at least $1-1/(2d)$.  So with positive probability all of them
do, and a spanning tree of each gives the pattern.  Applied to matching
graphs, with $\beta=\frac12$, $b=2s\le2n$ and $d=k$, it shows, as in the proof
of Theorem~\ref{thm:main} below, that almost every list of
$m^{\mathrm{HS}}_k(n)=k\,T(2n)$ directions separates multisets of size $n$.
Their proof needs only $(\log_2b)\bigl(1-\beta/(4\ln b)\bigr)^T\le1/(2d)$, and
$T(b)$ is a simple choice of $T$ that satisfies it with room to spare.  Using
that inequality as it stands, with nothing else changed, gives the smaller bound
\begin{equation}
\label{eq:mhse}
  m^{\mathrm{HSE}}_k(n)\;:=\;k\cdot\min\Bigl\{T\ge1:\
  \log_2(2n)\Bigl(1-\frac1{8\ln(2n)}\Bigr)^T\le\frac1{2k}\Bigr\},
\end{equation}
which we count as theirs (Remark~\ref{rem:sharp}).  Beyond that, the proof gives away constant factors in
four places.  Keeping them gives the following lemma, in which the four changes
are marked (i)--(iv).

\begin{lemma}[Hair and Sahai's lemma, tightened]
\label{lem:sharp}
Let $b\ge4$ be the number of vertices, $d\ge2$ be the number of spanning trees
required, let $\beta>0$ be the matching density factor, and let the integer
$T\ge1$ be the number of colours per tree, so that the total number of colours
is $m:=dT$.

Let $G$ be a loopless edge-coloured multigraph on $b$ vertices with at least $m$
colours, in which every colour class is a matching with at least $\beta b$
edges.  If
\begin{equation}
\label{eq:sharp}
  d\,\log_2b\,\prod_{j=0}^{T-1}\Bigl(1-\frac{\beta m}{2(m-j)\ln b}\Bigr)\;<\;1,
\end{equation}
then $G$ contains a pinning pattern with $d$ trees.
\end{lemma}

\begin{proof}
Discarding colours, we may assume that $G$ has exactly $m$ colours.  As in the
proof of Lemma~\ref{lem:hs}, choose a set $W$ of $w\ge2$ vertices that
maximises $e(U)/(|U|\ln|U|)$ over all sets $U$ of at least two vertices, and let
$\rho$ be the maximum.  Then $\rho\ge e(V(G))/(b\ln b)\ge\beta m/\ln b$.  Let
$H:=G[W]$.

\emph{Cuts.}  Let $A\subseteq W$ with $1\le a:=|A|\le w/2$.  By the choice of
$W$, $e(A)\le\rho\,a\ln a$ and $e(W\setminus A)\le\rho\,(w-a)\ln(w-a)$, so at
least $\rho\bigl(a\ln\frac wa+(w-a)\ln\frac w{w-a}\bigr)$ edges join $A$ to
$W\setminus A$.  A matching has at most $a$ of these edges, so at least
\[
  \rho\Bigl(\ln\frac wa+\frac{w-a}a\ln\frac w{w-a}\Bigr)\;\ge\;\rho\ln\frac{2w}a
\]
colours cross the cut.  For the inequality, put $t:=a/w\le\frac12$; then
$\frac{1-t}t\ln\frac1{1-t}=1-\sum_{i\ge1}\frac{t^i}{i(i+1)}$ decreases in $t$, so
it is at least its value $\ln2$ at $t=\frac12$.  (i)~Hair and Sahai drop the
second term, which leaves $\rho\ln(w/a)$.

\emph{One more colour.}  Let $c_1,\dots,c_m$ be the colours in a uniformly
random order, let $X_j$ be the number of components of
$H[\{c_1,\dots,c_j\}]$, and let $Z_j:=\ln X_j$, so that $Z_0=\ln w$.  Fix $j<m$
and $c_1,\dots,c_j$, and suppose that $x:=X_j\ge2$; probabilities and
expectations below are conditional on $c_1,\dots,c_j$.  Call a component
\emph{small} if it has at most $w/2$ vertices.  None of $c_1,\dots,c_j$ crosses
a component, so for a small component $C$ at least $\rho\ln(2w/|C|)$ of the
$m-j$ remaining colours cross it, and $c_{j+1}$ crosses it with probability at
least $\rho\ln(2w/|C|)/(m-j)$.  (ii)~Hair and Sahai replace $m-j$ by $m$ here.
So the number $Y$ of small components that $c_{j+1}$ crosses has
\[
  \mathbb E[Y]\;\ge\;\frac\rho{m-j}\sum_{C\text{ small}}\ln\frac{2w}{|C|}
  \;\ge\;\frac\rho{m-j}\,x\ln x .
\]
For the last step: if every component is small, Jensen's inequality gives
\[
  \sum_{C\text{ small}}\ln\frac{2w}{|C|}\;\ge\;x\ln(2x)\;\ge\;x\ln x .
\]
Otherwise exactly one of the components has more than $w/2$ vertices, so there
are $x-1$ small components with $s<w/2$ vertices in total, hence
\[
  \sum_{C\text{ small}}\ln\frac{2w}{|C|}\;\ge\;(x-1)\ln\frac{2w(x-1)}s
  \;\ge\;(x-1)\ln\bigl(4(x-1)\bigr)\;\ge\;x\ln x .
\]
The inequality $(x-1)\ln\bigl(4(x-1)\bigr)\ge x\ln x$ is justified, since the
two agree at $x=2$ and the derivative of the difference,
$\ln\bigl(4(x-1)/x\bigr)$, is positive on $x\ge2$.  (iii)~Hair and Sahai have
$\frac x2\ln x$ here.

Adding $c_{j+1}$ lowers the number of components by some $\Delta\ge Y/2$: in the
multigraph on the current components formed by the edges of colour $c_{j+1}$,
the components counted by $Y$ are not isolated, and $t$ nonisolated vertices lie
in at most $t/2$ nontrivial components.  Since
\[
  Z_j-Z_{j+1}\;=\;\ln\frac x{x-\Delta}\;\ge\;\frac\Delta x ,
\]
the preceding estimates give
\[
  \mathbb E[Z_{j+1}\mid c_1,\dots,c_j]\;\le\;
  Z_j-\frac{\mathbb E[\Delta\mid c_1,\dots,c_j]}x
  \;\le\;Z_j-\frac{\rho\ln x}{2(m-j)}
  \;\le\;\Bigl(1-\frac{\beta m}{2(m-j)\ln b}\Bigr)Z_j .
\]
This holds trivially when $x=1$.

\emph{$T$ colours.}  The factors in \eqref{eq:sharp} are positive, because
$\beta\le\frac12$ (a matching on $b$ vertices has at most $b/2$ edges),
$m/(m-j)<d/(d-1)\le2$ for $j<T$, and $\ln b\ge\ln4>1$.  Taking expectations and
iterating with $Z_0=\ln w$ yields
\[
  \mathbb E[Z_T]\;\le\;\ln w\prod_{j=0}^{T-1}\Bigl(1-\frac{\beta m}{2(m-j)\ln b}\Bigr).
\]
If $H[\{c_1,\dots,c_T\}]$ is disconnected, then $Z_T\ge\ln2$.  Thus Markov's
inequality gives
\[
  \mathbb P\bigl(H[\{c_1,\dots,c_T\}]\text{ is disconnected}\bigr)\;\le\;
  \frac{\mathbb E[Z_T]}{\ln2}\;\le\;
  \log_2b\prod_{j=0}^{T-1}\Bigl(1-\frac{\beta m}{2(m-j)\ln b}\Bigr),
\]
as $w\le b$.  Without (i)--(iii), the product is $(1-\beta/(4\ln b))^T$; Hair and Sahai bound it
by $e^{-\beta T/(4\ln b)}$ and then choose $T$ with room to spare.

\emph{$d$ trees.}  Split the order into $d$ consecutive blocks of $T$ colours.
Each block, in its order, is distributed like $c_1,\dots,c_T$, so by
\eqref{eq:sharp} and the union bound, with positive probability every block
connects $W$.  Spanning trees of the $d$ blocks then form the pattern.  (iv)~Hair
and Sahai ask each block to connect $W$ with probability at least $1-1/(2d)$; the
union bound needs only that the $d$ failure probabilities add up to less
than~$1$.
\end{proof}

\begin{proof}[Proof of Theorem~\ref{thm:main} for $m^\sharp_k(n)$ and $m^{\mathrm{CF}}_k(n)$]
Let $m\ge m^\sharp_k(n)=kT$, and let $D$ be a list of $m$ directions outside
the set of measure zero of Corollary~\ref{cor:generic-no-pattern}.  Suppose that
$D$ does not separate multisets of size $n$.  Then a smallest confusable pair
$X,Y$ has some size $s\le n$ and no point in common (Lemma~\ref{lem:pad}).

If $s=1$, its two points are joined by an edge of every colour, and any $k\le m$
colours give a pinning pattern on these two points, with one edge per tree.  If
$s\ge2$, its matching graph has $b=2s\le2n$ vertices and $m$ colours, and each
colour class is a perfect matching, with $b/2$ edges.  Lemma~\ref{lem:sharp}
applies with $\beta=\frac12$, $d=k$ and the $T$ of~\eqref{eq:msharp} (its $m$ is
$kT$, at most ours, and it may discard colours): the left side of \eqref{eq:sharp}
increases with $b$, since $\log_2b$ and each factor do, and it is below $1$ at
$b=2n$ by the choice of $T$.  In both cases the matching graph contains a
pinning pattern, which Corollary~\ref{cor:generic-no-pattern} rules out.

For the bound on $m^\sharp_k(n)$, let $m:=\bigl\lceil m^{\mathrm{CF}}_k(n)\bigr\rceil$,
and put $T:=\lfloor m/k\rfloor$, so that
$kT\ge m-k+1\ge(4k-2)\ln(2n)\ln\bigl(k\log_2(2n)\bigr)+2$.  It is enough to show
that this $T$ satisfies the condition that defines $m^\sharp_k(n)$, for then
$m^\sharp_k(n)\le kT\le m$.  By
$1-u\le e^{-u}$, it does if
\[
  S\;:=\;\sum_{j=0}^{T-1}\frac{kT}{kT-j}\;>\;4\ln(2n)\,\ln\bigl(k\log_2(2n)\bigr).
\]
Now
\[
  S\;\ge\;kT\ln\frac{kT+1}{(k-1)T+1}\;\ge\;kT\ln\frac k{k-1}-\frac1{k-1}
  \;\ge\;\frac{2kT}{2k-1}-\frac1{k-1} .
\]
The first step is $\sum_{i=p+1}^q\frac1i\ge\ln\frac{q+1}{p+1}$.  For the second,
$\frac{kT+1}{(k-1)T+1}=\frac k{k-1}(1-u)$ with $u=\frac1{k(k-1)T+k}$, and
$\ln(1-u)\ge-\frac u{1-u}$.  The third is $\ln(1+y)\ge\frac{2y}{2+y}$ with
$y=\frac1{k-1}$.  With the bound on $kT$, $S$ is at least
$4\ln(2n)\ln\bigl(k\log_2(2n)\bigr)+\frac4{2k-1}-\frac1{k-1}$, and
$\frac4{2k-1}>\frac1{k-1}$ for $k\ge2$.

Finally, $m^{\mathrm{CF}}_k(n)$ is $O\bigl(k\log n\,(\log\log n+\log k)\bigr)$, and
since almost every list of that many directions separates, certainly some list
does.
\end{proof}

\begin{remark}[Where the saving comes from]
\label{rem:sharp}
For $k=2$ and $n=10^3$, Hair and Sahai's lemmas as their proof writes them give
$m^{\mathrm{HS}}_2(n)=558$, and their final inequality evaluated exactly gives
$m^{\mathrm{HSE}}_2(n)=458$ (equation~\eqref{eq:mhse}).  We count both as theirs, carried over to sorted projections.  Our
change~(iv) then gives $374$, (ii) gives $270$, and (i) with (iii) gives
$m^\sharp_2(n)=134$, so our changes save a factor of about $3.4$.  Letting some
blocks fail (Lemma~\ref{lem:blocks-fail}) then gives
$m^\flat_{2,\mathrm{bal}}(n)=133$.  The closed form (equation~\eqref{eq:mcf}) gives $\bigl\lceil m^{\mathrm{CF}}_2(n)\bigr\rceil=144$.
\end{remark}

\paragraph{Allowing some blocks to fail: $m^\flat_{k,\mathrm{eq}}(n)$ and $m^\flat_{k,\mathrm{bal}}(n)$.}
In Lemma~\ref{lem:sharp}, the $m=dT$ colours are partitioned into $d$ blocks of
$T$ colours each, and the union bound requires every one of the $d$ blocks to
connect $W$.  More generally, we can partition $m$ colours into $q\ge d$
disjoint blocks of lengths $T_1,\dots,T_q$ and ask only that at least $d$ of
them connect $W$, allowing up to $q-d$ blocks to fail.  This gives the stronger
bounds $m^\flat_{k,\mathrm{eq}}(n)$ and $m^\flat_{k,\mathrm{bal}}(n)$ of
Theorem~\ref{thm:main}.

\begin{lemma}[Allowing blocks to fail]
\label{lem:blocks-fail}
Let $b\ge4$, $d\ge2$, and $0<\beta\le\frac12$.  Let $q\ge d$ and let
$T_1,\dots,T_q\ge1$ be integers with $\sum_{l=1}^q T_l=m$ and
$T_l\le\lceil m/2\rceil$ for each $l$.  Let $G$ be a loopless edge-coloured
multigraph on $b$ vertices with at least $m$ colours, in which every colour
class is a matching with at least $\beta b$ edges.  If
\begin{equation}
\label{eq:blocks-fail}
  \sum_{l=1}^q \log_2b\,\prod_{j=0}^{T_l-1}\Bigl(1-\frac{\beta m}{2(m-j)\ln b}\Bigr)\;<\;q-d+1,
\end{equation}
then $G$ contains a pinning pattern with $d$ trees.
\end{lemma}

\begin{proof}
Discarding excess colours, we may assume that $G$ has $m$ colours, and choose
$W$ of size $w\ge2$ as in the proof of Lemma~\ref{lem:sharp} before ordering
the colours.  Split a uniformly random ordering of the $m$ colours into $q$
consecutive blocks of lengths $T_1,\dots,T_q$.  For each $l$, since
$j<T_l\le\lceil m/2\rceil$, we have $m/(m-j)<2$ and
$\beta m/\bigl(2(m-j)\ln b\bigr)<1/(2\ln 4)<1$, so every factor in
\eqref{eq:blocks-fail} is positive.  Each block of $T_l$ colours has the same
marginal distribution as the first $T_l$ colours $c_1,\dots,c_{T_l}$,
regardless of its place in the ordering, so by the single-block estimate in the
proof of Lemma~\ref{lem:sharp} the $l$-th block fails to connect $W$ with
probability at most the $l$-th term of \eqref{eq:blocks-fail}.

Let $F$ be the number of blocks that fail to connect $W$.  By linearity of
expectation and \eqref{eq:blocks-fail}, $\mathbb E[F]<q-d+1$.  Since $F$ takes
integer values, with positive probability $F\le q-d$, so at least $q-F\ge d$
blocks connect $W$.  Spanning trees of $d$ of those blocks form a pinning
pattern with $d$ trees.
\end{proof}

\begin{proof}[Proof of Theorem~\ref{thm:main} for $m^\flat_{k,\mathrm{eq}}(n)$ and $m^\flat_{k,\mathrm{bal}}(n)$]
In \eqref{eq:mflateq}, setting $q=k$ recovers the definition \eqref{eq:msharp}
of $m^\sharp_k(n)$, so $m^\flat_{k,\mathrm{eq}}(n)\le m^\sharp_k(n)$.  In
\eqref{eq:mflatbal}, restricting $m$ to multiples $qT$ of $q$ (so $r=0$ and
$B_{2n}(qT,q)=q\,P_{2n}(qT,T)$) recovers \eqref{eq:mflateq}, so
$m^\flat_{k,\mathrm{bal}}(n)\le m^\flat_{k,\mathrm{eq}}(n)$.

Now let $m\ge m_0:=m^\flat_{k,\mathrm{bal}}(n)$, and choose $q\in\{k,\dots,m_0\}$
with $B_{2n}(m_0,q)<q-k+1$.  Write $m_0=qt+r$ with $0\le r<q$, and take $q-r$
blocks of length $t\ge1$ and $r$ blocks of length $t+1$; since $q\ge k\ge2$,
every block has length at most $\lceil m_0/q\rceil\le\lceil m_0/2\rceil$.  As
above, for any confusable pair of size $s\le n$ with no point in common, if
$s=1$ then any $k\le m$ colours give a pinning pattern, and if $s\ge2$ then
$b=2s\le2n$ and $P_b(m_0,\cdot)\le P_{2n}(m_0,\cdot)$, so
Lemma~\ref{lem:blocks-fail} (with $\beta=\frac12$, $d=k$, after discarding
$m-m_0$ colours) supplies a pinning pattern on $W$, which
Corollary~\ref{cor:generic-no-pattern} rules out for generic lists.
\end{proof}

Tables~\ref{tab:numbers} and~\ref{tab:flat-comparison} in the appendix compare
$m^\sharp_k(n)$, $m^\flat_{k,\mathrm{eq}}(n)$ and $m^\flat_{k,\mathrm{bal}}(n)$
numerically across $n$ and $k$.

\begin{remark}[High-dimensional bound]
\label{rem:no-logk}
By $1-u\le e^{-u}$ and $m/(m-j)\ge1$, we have
$P_N(qT,T)\le(\log_2N)\,e^{-T/(4\ln N)}$, so a block of length
$T=\bigl\lceil4\ln N\,\ln\bigl(q\log_2N/(q-k+1)\bigr)\bigr\rceil$ has failure
probability at most $(q-k+1)/q$.  Taking $q=2k$, $N=2n$ and block length
$T=\bigl\lceil4\ln(2n)\,\ln\bigl(2\log_2(2n)\bigr)\bigr\rceil$ gives
$P_{2n}(2kT,T)\le\frac12$, so $2k\,P_{2n}(2kT,T)\le k<k+1=q-k+1$.  This yields
the closed-form \emph{high-dimensional bound}
\begin{equation}
\label{eq:mhd}
  m^{\mathrm{HD}}_k(n)\;:=\;2k\bigl\lceil4\ln(2n)\,\ln\bigl(2\log_2(2n)\bigr)\bigr\rceil,
\end{equation}
for which $m^\flat_{k,\mathrm{bal}}(n)\le m^\flat_{k,\mathrm{eq}}(n)\le m^{\mathrm{HD}}_k(n)$.
Unlike $m^{\mathrm{CF}}_k(n)$ and $m^\sharp_k(n)$, this bound has no $\log k$
factor---it gives $m_k(n)=O(k\log n\,\log\log n)$ uniformly in $k$ and $n$---and
so scales better when the dimension $k$ is large and $n$ is small.
\end{remark}

\section{Almost every list behaves alike}
\label{sec:generic}

The proof of Theorem~\ref{thm:zero-one} applies Lemma~\ref{lem:criterion}
with the coordinates of the directions as unknowns.  It uses the notation set
up just before that lemma: a list $\rho=(\rho_1,\dots,\rho_m)$ of permutations
of $\{1,\dots,n\}$ gives the linear equations
$x_a\cdot d_j=y_{\rho_j(a)}\cdot d_j$ ($a=1,\dots,n$, $j=1,\dots,m$), whose
matrix is $M_\rho(D)$ and whose solutions form $K_\rho(D)$; and $Z_\sigma$ is
the subspace on which $x_a=y_{\sigma(a)}$ for every $a$.

\begin{proof}[Proof of Theorem~\ref{thm:zero-one}]
Fix $\rho$.  The entries of the matrix $M_\rho(D)$ are $0$ or $\pm$ a
coordinate of a direction.  Let $r$ be the largest
rank of $M_\rho(D)$ over all lists $D$.  A list where it is attained is a
witness (Section~\ref{sec:prelim}) that some $r\times r$ minor is a polynomial
in $D$ that is not identically zero, so the rank is $r$ on a set $U_\rho$ of
full measure.  For $D$ in $U_\rho$, $K_\rho(D)\subseteq Z_\sigma$ exactly when
adding the equations of $Z_\sigma$ to $M_\rho(D)$ does not raise its rank
(Lemma~\ref{lem:criterion}), that
is, when every $(r+1)\times(r+1)$ minor of the enlarged matrix vanishes.  If, for
some $\sigma$, these minors all vanish identically, then $K_\rho(D)$ lies in
$Z_\sigma$ for every $D$ in $U_\rho$.  Otherwise, for each $\sigma$ one of them
is a nonzero polynomial, and off the union of their zero sets $K_\rho(D)$ lies
in no $Z_\sigma$.  So the lists for which $K_\rho(D)$ lies in some $Z_\sigma$
form a set of full measure or a set of measure zero.  By
Lemma~\ref{lem:criterion}, the lists
that separate multisets of size $n$ are the intersection of these sets over the finitely many
$\rho$, so they too form a set of full measure or of measure zero.

All these polynomials have integer coefficients.  At a list whose $mk$
coordinates are algebraically independent, such a polynomial vanishes only if
it is identically zero, so such a list is a witness for all of them at once.  So such a list lies in every $U_\rho$, and, for each
$\rho$, $K_\rho(D)$ lies in some $Z_\sigma$ exactly when this holds for almost
every list.  So it separates multisets of size $n$ exactly when almost every
list does.
\end{proof}

\begin{corollary}[One random sequence serves every $n$]
\label{cor:nested}
Let $d_1,d_2,\dots$ be random directions such that, for every $m$, the list
$(d_1,\dots,d_m)$ is drawn as in (ii) or (iii) of Table~\ref{tab:random}; for
example, independent directions with independent standard normal coordinates,
or uniform on the unit sphere.  With probability one, for every $n\ge2$ and every
$m\ge m^\flat_{k,\mathrm{bal}}(n)$ (and in particular every $m\ge m^\sharp_k(n)$),
the list $(d_1,\dots,d_m)$ separates multisets of size $n$.
\end{corollary}

\begin{proof}
For each of the countably many pairs $(m,n)$, the first $m$ directions fail
with probability zero, by Theorem~\ref{thm:main}, which holds for generic
lists, and Lemma~\ref{lem:random}(c).  A countable union of events of probability zero has probability zero.
\end{proof}

It can be convenient to fix some of the directions, for example to make the
first $k$ of them the coordinate axes, so that $\Phi_D(X)$ begins with the
sorted coordinates of the points, as in Remark~\ref{rem:small} and
\cite[Section~4]{DWTLB}.  Fixing them in this way does not change which
results hold for almost every list, as the next corollary shows.

\begin{corollary}[The axes may be among the directions]
\label{cor:axes}
Fix $k$, $m\ge k$ and $n$.  The list $(e_1,\dots,e_k,d_{k+1},\dots,d_m)$
separates multisets of size $n$ for almost every choice of
$d_{k+1},\dots,d_m$ (Definition~\ref{def:ae}) if and only if almost every list
of $m$ directions does.  In particular it does so when
$m\ge m^\flat_{k,\mathrm{bal}}(n)$ (or $m\ge m^\sharp_k(n)$).
\end{corollary}

\begin{proof}
For an invertible matrix $A$, a list $D$ separates multisets of size $n$ exactly
when $AD=(Ad_1,\dots,Ad_m)$ does, because the projection of $X$ on $Ad$ is the
projection of $A^{\mathsf T}X$ on $d$.  The map
$(A,d_{k+1},\dots,d_m)\mapsto(Ae_1,\dots,Ae_k,Ad_{k+1},\dots,Ad_m)$ is a
diffeomorphism from $GL_k(\R)\times(\R^k)^{m-k}$ onto the set of lists whose
first $k$ directions are linearly independent, which has full measure.  The
preimage of the set $B$ of lists that fail is $GL_k(\R)\times B'$, where $B'$ is
the set of $(d_{k+1},\dots,d_m)$ for which $(e_1,\dots,e_k,d_{k+1},\dots,d_m)$
fails.  By Fubini's theorem $B'$ has measure zero exactly when $B$ does.
\end{proof}

The two corollaries combine: by the same countable union, with probability one
the axes followed by one random sequence $d_{k+1},d_{k+2},\dots$ serve every
$n$.

\begin{remark}[Compression]
\label{rem:compress}
The map $\Phi_D$ has $mn$ outputs, whereas a multiset of $n$ points in $\R^k$
has only $nk$ coordinates.  Applying a further linear map after $\Phi_D$ can
reduce the outputs to $(2n-1)k$, whatever $m$ is, and almost any such map will
do: if $D$ separates multisets of size $n$, then for
almost every linear map $L:\R^{mn}\to\R^{(2n-1)k}$ (Definition~\ref{def:ae})
the map
$X\mapsto L(\Phi_D(X))$ is still injective on multisets of size $n$.  This is
the argument of \cite[Theorem~10]{DWTLB}, which is stated for $(k-1)n+1$
directions in general position but uses that hypothesis only to know that
$\Phi_D$ is injective.  In outline: for fixed orders of the points of $X$ and
of $Y$ in each sorted list, $\Phi_D(X)-\Phi_D(Y)$ is a linear function of the
$2n$ points that vanishes when they are all equal, so it takes values in one of
finitely many subspaces of dimension at most $(2n-1)k$.  Almost every $L$ is
injective on each of them, and then $L\Phi_D(X)=L\Phi_D(Y)$ forces
$\Phi_D(X)=\Phi_D(Y)$, hence $X=Y$.  This is the finite witness theorem
(Section~\ref{sec:prelim}) in its linear form, with the rows of $L$ as the
witnesses; \cite{DWTLB} prove their Theorem~10 with it.  The compressed map is then bi-Lipschitz too
\cite[Theorem~1]{BTW24}, and \cite[Theorem~1.2]{BHS} already has a version with
$2nk$ outputs.
\end{remark}

\section{Open questions}
\label{sec:open}

For fixed $k$, Theorem~\ref{thm:main} and the lower bounds differ by a factor
of order $\log\log n$; for the $n$ of Table~\ref{tab:numbers} the factor is
between $7$ and $16$.

\begin{itemize}[leftmargin=*,itemsep=2pt]
\item Is $m_k(n)=O(k\log n)$?  For fixed $k$, this asks whether the largest $n$
  that some list of $m$ directions separates grows exponentially with $m$.  The
  bounds of \cite{MPS} leave this open even for $k=2$.
\item How many colours force a pinning pattern?  Let $g_k(N)$ be the largest
  number of colours of an edge-coloured bipartite multigraph on $N$ vertices,
  with every colour class a perfect matching, that contains no pinning pattern.
  By the proof of Theorem~\ref{thm:main}, almost every list of $m$ directions
  separates multisets of size $n$ once $m>g_k(2s)$ for every $s\le n$, so
  $g_k(N)=O(k\log N)$ would answer the first question, even for almost every
  list.  For $k=2$, \cite[Section~3.4]{MPS} show that, under the cut condition
  of their Lemma~3.7, choosing the colours of the trees at random can fail even
  with about $\log n\,\log\log n/\log\log\log n$ colours, although the graphs of
  their example do contain the trees.  An exhaustive search, in the accompanying code, gives
  $g_2(10)=4$, whereas for $k=2$ and $N=10$ Lemma~\ref{lem:sharp} needs $24$
  colours.
\item Every list of five directions in the plane fails on multisets of size
  six (Lemma~\ref{lem:five}), and this gives the rate $12^{1/5}\approx1.644$
  of Theorem~\ref{thm:plane}.  Does every list of seven directions fail on
  multisets of size sixteen?  By Lemma~\ref{lem:product}, that would give a
  smaller rate.  Six directions cannot give one (Remark~\ref{rem:six-seven}),
  and almost every list of seven separates multisets of size $14$, so only the
  sizes $15$ and $16$ remain.  The opposite may be true: the lower bound of
  Theorem~\ref{thm:plane} equals $m_2(n)$ for every $n\le14$
  (Remarks~\ref{rem:small} and~\ref{rem:six-seven}), and if it did so for
  every $n$, then $m_2(n)=5\log_{12}n+O(1)$.  So: for every $m=5t+r$ with
  $0\le r\le4$, and $c_r$ as in Theorem~\ref{thm:plane}, does some list of $m$
  directions in the plane separate multisets of size $\frac12c_r12^t-1$?
\item Is $m_k(3)=2k$ for every $k$?  It holds for $k\le4$
  (Theorem~\ref{thm:three} and Remark~\ref{rem:small}).
\item Is there ever a list that separates multisets of size $n$ although almost
  every list of the same length does not?  By Theorem~\ref{thm:zero-one}, its
  coordinates would be algebraically dependent.  For given $k$, $m$ and $n$,
  whether almost every list separates can be decided by running the test of
  Lemma~\ref{lem:criterion} with indeterminates as coordinates.  In the plane
  there is no such list for $n\le14$ (Remark~\ref{rem:six-seven}).
\item Is there a simple explicit list of $O(k\,\mathrm{polylog}\,n)$ directions
  that separates multisets of size $n$?  By Theorems~\ref{thm:main}
  and~\ref{thm:zero-one}, a list of $m\ge m^\sharp_k(n)$
  directions with algebraically independent coordinates does, but, as
  \cite[Section~5]{MPS} remark for $k=2$, such a list is hardly simple.  Their
  work was motivated by a question of Holub: how large must $m$ be for $n$
  functions $x\mapsto a_ix+b_i$ to be determined by the multisets of their
  values at $x=1,\dots,m$?  In our terms this concerns the directions
  $(1,1),(1,2),\dots,(1,m)$ in the plane, and for them neither their methods
  nor ours improve on the linear bound.
\item How large is the distortion of $\Phi_D$, the ratio of its upper and lower
  Lipschitz constants, when $D$ has only $O(k\,\mathrm{polylog}\,n)$
  directions?  With about $kn^2\log n$ directions drawn uniformly from the unit
  sphere it is $O(n^2)$ up to
  logarithmic factors, and for every list it is at least of order $\sqrt n$
  \cite[Section~3]{DWTLB}.  Can the compression of Remark~\ref{rem:compress}
  use a linear map with few nonzero entries in each row?  A dense one costs
  more to apply than $\Phi_D$ itself.
\end{itemize}

\section*{Acknowledgements}
CGL wishes to thank Peterhouse, not only for providing financial support for
the work undertaken, but also for providing an environment that has fostered
many positive discussions that have influenced this work.  Notably,
contributions from Patrick Kennedy-Hunt (on simplicial encoders) and Ian Wright
(on polynomial encoders) were instrumental in shaping our interest in
bi-Lipschitz encoding methods.

CGL also wishes to acknowledge Tom Ruane (then a Durham undergraduate) for
productive correspondence over the summer of 2025.

\section*{Statement on AI use}
This work has a long gestation beginning in 2019.\footnote{Arguably this work's
beginnings pre-date CGL's work in 2020 on permutation invariant encodings
\cite{GHL21,LHG22}.  Observations of deficiencies in those encoding strategies
led to CGL's increased interest in projection sort encoders in late 2023, and
the subsequent creation of \url{https://github.com/kesterlester/Echinocoder/} in
February 2024.  A fruitful association with a summer student (Tom Ruane)
followed in 2025, and RS was recruited as an MPhil student in early June 2026
for an October 2026 start.}  Prior to 2025, no LLM made any significant
contribution to this work at all.  Between 2025 and 21~October 2026, LLMs were
used as tools to assist with simple code refactoring tasks under complete
direction.  From 21~October 2026 until submission to the arXiv, substantial use
was made of Anthropic's Claude Opus~5.5 to assist in the drafting of this paper
and to provide both criticism of and refinements to arguments and methods it
contains---although the authors have checked every statement and take
responsibility for the work.  This sudden switch was a direct consequence of
reading the statement regarding AI use in \cite{HS26}, which itself followed
our discovery that results in \cite{HS26} had rendered obsolete many of the
things we \emph{had} hoped to write about as extensions to \cite{DWTLB}
and~\cite{MPS}.  Separately, some OpenAI models were used to assist with the
discovery of relevant literature.

\appendix
\section{Conditions on the directions}
\label{app:conditions}

\begin{figure}[ht]
\centering
\begin{tikzpicture}[font=\footnotesize,
    sep/.style={fill=batB!22,draw=batB!75!black,thick},
    oval/.style={draw=black,thin}]
  \filldraw[sep] (4.2,1.7) ellipse (2.4 and 1.1);
  \filldraw[sep] (7.6+3.85,1.93) ellipse (3.0 and 1.58);
  \foreach \xs in {0,7.6} {
    \begin{scope}[xshift=\xs cm]
      \draw[gray] (0,0) rectangle (7,4.6);
      \node[anchor=north west,gray] at (0,4.6) {lists of $m$ directions};
      \draw[oval] (3.3,2.2) ellipse (2.9 and 1.85);
      \node at (3.3,3.7) {non-parallel};
      \draw[oval] (3.3,2.0) ellipse (2.2 and 1.3);
      \node at (3.3,2.98) {general position};
      \draw[oval,dashed] (3.3,1.75) ellipse (1.1 and 0.6);
      \node[align=center,font=\scriptsize] at (3.3,1.75)
        {algebraically\\[-1pt] independent};
    \end{scope}}
  \node[text=batB!60!black] at (5.0,0.95) {$S$};
  \node[text=batB!60!black] at (7.6+5.2,1.15) {$S$};
  \node at (3.5,-0.45) {(a) $ m^\sharp_k(n)\le
    m\le(k-1)n$};
  \node at (7.6+3.5,-0.45) {(b) $m\ge(k-1)n+1$};
\end{tikzpicture}
\caption{The conditions of this appendix on a list of $m\ge k$ directions in
$\R^k$, drawn as sets for $k\ge3$ (for $k=2$ the two outer ovals coincide),
and the set $S$ of lists that separate multisets of size $n$, shaded.
(a)~When $ m^\sharp_k(n)\le m\le(k-1)n$, $S$ contains
every list with algebraically independent coordinates
(Theorem~\ref{thm:main} and Lemma~\ref{lem:generic}), but not every list in
general position (Remark~\ref{rem:gp-sharp}).  (b)~When $m\ge(k-1)n+1$, $S$
contains every list in general position (Lemma~\ref{lem:linear}).  When a
shorter list separates, $S$ also contains lists with two parallel directions
(repeat a direction of the shorter list), so $S$ reaches outside the outer
oval.  Areas mean nothing: every set drawn has full measure, so the sets differ
only by sets of measure zero.}
\label{fig:conditions}
\end{figure}
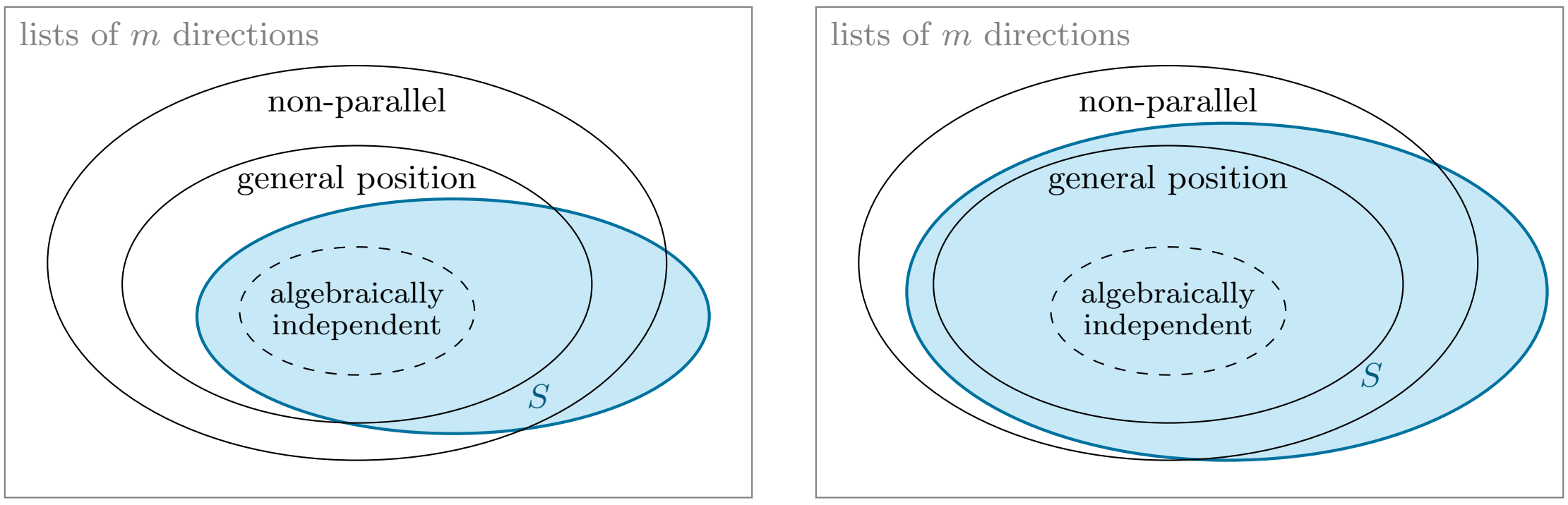

The results of this paper put conditions of three kinds on a list
$D=(d_1,\dots,d_m)$ of directions in $\R^k$: properties of the one list $D$,
statements about almost all lists, and ways of drawing $D$ at random.  This
appendix defines them, compares them (Figure~\ref{fig:conditions}), and
records which result uses which (Table~\ref{tab:uses}).  Directions are nonzero
vectors, and every property in this paper depends only on the lines through
them, since $\pi_{\lambda d}X$ is $\pi_dX$ scaled by $\lambda$; so each
condition below is unchanged when a direction is replaced by a nonzero
multiple of itself.  Two directions are \emph{parallel} if one is a multiple of
the other.

\paragraph{General position.}

\begin{definition}[General position]
\label{def:gp}
A list $d_1,\dots,d_m$ of $m\ge k$ directions in $\R^k$ is in \emph{general
position} if every $k$ of $d_1,\dots,d_m$ are linearly independent.
\end{definition}

\begin{lemma}
\label{lem:gp}
For a list $D=(d_1,\dots,d_m)$ of $m\ge k$ directions in $\R^k$ the following
are equivalent:
\begin{enumerate}[label=(\alph*)]
\item $D$ is in general position;
\item every $j\le k$ of $d_1,\dots,d_m$ are linearly independent;
\item every nonzero vector of $\R^k$ is orthogonal to at most $k-1$ of
  $d_1,\dots,d_m$.
\end{enumerate}
In particular the directions of a list in general position are pairwise
non-parallel, and for $k=2$ the converse holds.
\end{lemma}

\begin{proof}
Since $m\ge k$, any $j\le k$ of the directions are among some $k$ of them, and
part of an independent set is independent; so (a) implies (b), and (a) is the
case $j=k$ of (b).  And $k$ vectors of $\R^k$ are linearly dependent exactly
when they lie in a hyperplane through the origin, that is, when some nonzero
vector is orthogonal to all of them; so (a) and (c) are equivalent.  The last
statement is (b) with $j=2$, which for $k=2$ is (a).
\end{proof}

The condition $m\ge k$ matters.  For fewer than $k$ directions, (a) and (c)
would hold for every list, for instance for $(1,0,0),(2,0,0)$ in $\R^3$, while
(b) would ask that all of them be linearly independent; we never use the term
for fewer than $k$ directions.  In terms of the lines through the directions,
(b) says that the $m$ corresponding points of the projective space
$\mathbb{RP}^{k-1}$ are in general position in the usual sense: for
$2\le g\le k$, no $g$ of them lie in a projective subspace of dimension $g-2$.
This is not general position of $d_1,\dots,d_m$ as points of $\R^k$:
$(1,0,0)$ and $(2,0,0)$ are different points but the same direction.
Definition~\ref{def:gp} is the hypothesis of \cite[Proposition~48]{HSW}, which
concerns $(k-1)N+1\ge k$ unit vectors, and it is the \emph{full spark}
condition of \cite{DWTLB}, which is defined for matrices with at least as many
columns as rows.  Their proof of \cite[Theorem~4]{DWTLB}, like that of
Lemma~\ref{lem:linear}, uses it in the form~(c).

Whether a given list is in general position can be decided by computing its
$\binom mk$ determinants of $k$ directions, in exact arithmetic when the
directions are rational.  Explicit lists in general position are easy to write
down: the directions $(1,t_j,t_j^2,\dots,t_j^{k-1})$ with distinct numbers
$t_1,\dots,t_m$, since any $k$ of them form a Vandermonde matrix.  For $k=2$
and $t_j=j$ these are the directions of Holub's question
(Section~\ref{sec:open}).  General position is exactly what the linear bound
needs, and that bound cannot be lowered for it:

\begin{remark}[The linear bound is tight for general position]
\label{rem:gp-sharp}
For all $k\ge2$ and $n\ge2$ there is a list of $(k-1)n$ directions in general
position that does not separate multisets of size $n$.  For $k=2$, take the $n$
directions at angles $(2l+1)\pi/(2n)$, $l=0,\dots,n-1$, and let $X$ and $Y$ be
the vertices of a regular $n$-gon centred at the origin and of its rotation by
$\pi/n$.  The reflection in the line of each direction $d$ maps $X$ onto $Y$
and preserves $x\cdot d$, so $X$ and $Y$ have the same projection on $d$; this
is the example of \cite[Fig.~1]{BL90} and \cite[Section~1]{MPS}.  For $k\ge3$, write
$\R^k=\R^2\times\R^{k-2}$, let $u_1,\dots,u_n$ be these planar directions, and
put $X'=X\times\{0\}$ and $Y'=Y\times\{0\}$.  Then $X'$ and $Y'$ have the same
projection on every direction in each of the hyperplanes
$H_i:=\{(\lambda u_i,w):\lambda\in\R,\ w\in\R^{k-2}\}$.  Choose $k-1$
directions in each $H_i$.  Any $k$ of the $(k-1)n$ directions include two from
different hyperplanes $H_a$ and $H_b$, and their determinant, as a polynomial
in the parameters $\lambda$ and $w$ of the chosen directions, is not
identically zero, as a witness shows: it is $\pm\det(u_a,u_b)\neq0$ when those two are $(u_a,0)$
and $(u_b,0)$ and the others are $(0,e_1),\dots,(0,e_{k-2})$.  So the
parameters for which one of these determinants vanishes form a set of measure
zero, and any other choice gives a list in general position.  Dropping
directions keeps a list in general position, if at least $k$ remain, and keeps
it failing, so the same holds for every length $m$ with $k\le m\le(k-1)n$.
\end{remark}

\paragraph{What ``almost every'' and ``generic'' mean.}

\begin{definition}[Almost every; generic]
\label{def:ae}
Fix $k$ and $m$.  A property of lists of $m$ directions in $\R^k$ holds for
\emph{almost every} list if the lists that fail it form a set of Lebesgue
measure zero in $\R^{mk}$.  It holds for \emph{generic} lists if the lists that
fail it lie in the zero set of a polynomial in the $mk$ coordinates, with
integer coefficients, that is not identically zero.  When some directions are
fixed, as in Corollary~\ref{cor:axes}, or confined to given subspaces, as in
Remarks~\ref{rem:small} and~\ref{rem:gp-sharp}, ``almost every'' refers to
Lebesgue measure on the parameters that remain free; for the linear maps of
Remark~\ref{rem:compress}, to Lebesgue measure on their entries.
\end{definition}

\begin{definition}[Algebraically independent]
\label{def:ai}
Real numbers are \emph{algebraically independent} if no polynomial with integer
coefficients that is not identically zero vanishes at them.
\end{definition}

\begin{lemma}
\label{lem:generic}
If a property holds for generic lists of $m$ directions, then
\begin{enumerate}[label=(\alph*)]
\item it holds for almost every list;
\item it holds on a dense open set of lists;
\item it holds for every list whose $mk$ coordinates are algebraically
  independent.
\end{enumerate}
\end{lemma}

\begin{proof}
Let $P$ be the polynomial of Definition~\ref{def:ae}.  (a)~The zero set of a
polynomial that is not identically zero has measure zero, by induction on the
number of variables.  Write $P$ as a polynomial in the last variable whose
coefficients are polynomials in the others; one of them, $c$ say, is not
identically zero.  By induction $c$ vanishes only on a set of measure zero, and
wherever $c\neq0$, $P$ has finitely many roots in the last variable.  So by
Fubini's theorem the zero set of $P$ has measure zero.  (b)~The zero set of $P$
is closed, and by~(a) it contains no open set.  (c)~$P$ does not vanish at such
a list.
\end{proof}

None of these implications can be reversed, and (b) and (c) describe
different sets.  The property that the first coordinate of $d_1$ is irrational
holds for almost every list, but on no open set.  The lists with algebraically
independent coordinates have full measure, but they contain no open set, since
lists with rational coordinates are dense.  Corollary~\ref{cor:generic-no-pattern},
and with it Theorem~\ref{thm:main}, is proved in the generic form, and so is
Theorem~\ref{thm:zero-one}: for given $k$, $m$ and $n$, either generic lists
separate multisets of size $n$, or the lists that separate lie in the zero set
of a polynomial with integer coefficients that is not identically zero.  The
second statement of Theorem~\ref{thm:zero-one} rests on this form, and it is
how \cite[Theorem~1.1]{MPS} state their bound.

\paragraph{Random directions.}  Table~\ref{tab:random} lists four ways of
drawing a list at random.

\begin{table}[ht]
\centering
\small
\begin{tabular}{@{}lp{7.1cm}ccc@{}}
  \toprule
  & The list $D$ of $m\ge k$ directions is drawn & general & generic &
    almost every\\
  & & position & lists & list\\
  \midrule
  (i) & with independent directions, each from a distribution that gives no
    mass to hyperplanes through the origin & yes & if $k=2$ & ---\\
  (ii) & with a density on $\R^{mk}$, or with the unit vectors
    $d_j/\lvert d_j\rvert$ having a density on the product of $m$ unit
    spheres & yes & yes & yes\\
  (iii) & with independent coordinates, none of which takes any single value
    with positive probability & yes & yes & ---\\
  (iv) & with rounded coordinates, or from a finite set of lists & --- & --- &
    ---\\
  \bottomrule
\end{tabular}
\caption{Ways of drawing random directions, and what holds with probability
one (Lemma~\ref{lem:random}).  ``Yes'' under ``generic lists'' or ``almost every
list'' means that every property of the paper that holds for generic lists, or
for almost every list, then holds with probability one.  A dash means that
nothing is guaranteed.  Independent standard normal coordinates satisfy
(i)--(iii); independent directions uniform on the unit sphere satisfy (i)
and~(ii); independent directions uniform on a small circle of the unit sphere
of $\R^3$ satisfy only~(i).}
\label{tab:random}
\end{table}

\begin{lemma}[Random directions]
\label{lem:random}
Let $Q$ be a property of lists of $m$ directions that depends only on the lines
through them, as every property in this paper does, and let $D$ be drawn in one
of the ways (i)--(iii) of Table~\ref{tab:random}.
\begin{enumerate}[label=(\alph*)]
\item If $m\ge k$, then $D$ is in general position with probability one.
\item If $Q$ holds for almost every list and $D$ is drawn as in (ii), then $D$
  has $Q$ with probability one.
\item If $Q$ holds for generic lists and $D$ is drawn as in (ii) or (iii), then
  $D$ has $Q$ with probability one.  The same holds if $D$ is drawn as in (i)
  and $k=2$, provided that the polynomial of Definition~\ref{def:ae} is
  homogeneous in each direction separately.
\end{enumerate}
\end{lemma}

\begin{proof}
(b)~The lists that fail $Q$ form a set $B$ of measure zero that contains,
with each list, every list obtained by rescaling its directions.  A density on
$\R^{mk}$ gives $B$ probability zero.  In polar coordinates
$d_j=r_js_j$, Lebesgue measure on $(\R^k\setminus\{0\})^m$ is the product of
the measures $r_j^{k-1}\,dr_j\,d\sigma(s_j)$, where $\sigma$ is surface measure
on the unit sphere, so $B$ meets the product of the spheres in a set of surface
measure zero, and a density there gives it probability zero too.

(c)~Under (ii) this follows from (b) and Lemma~\ref{lem:generic}(a).  Under
(iii), it follows from the induction in the proof of
Lemma~\ref{lem:generic}(a), with ``measure zero'' replaced by ``probability
zero'': a polynomial in one variable that is not identically zero has finitely
many roots, and the last coordinate takes each of them with probability zero.
Under (i) with $k=2$, the same induction runs over the directions instead of
the coordinates: for given values of the other directions, a polynomial that is
homogeneous in $d_j\in\R^2$ and not identically zero vanishes only on finitely
many lines through the origin, and $d_j$ misses them with probability one.

(a)~Under (i), let $j<k$ of the directions be linearly independent.  They lie
in some hyperplane through the origin, and a further direction, drawn
independently, lies in it with probability zero.  So any $k$ given directions
are independent with probability one, and there are finitely many choices of
$k$ of them.  Under (ii) and (iii), general position holds for generic lists,
since a list that is not in general position is a zero of the product of all
the determinants of $k$ of its directions; so (a) follows from~(c).
\end{proof}

All the polynomials of Sections~\ref{sec:rigidity} and~\ref{sec:generic} are
homogeneous in each direction, because each row of the matrices there is
linear in a single direction, so Lemma~\ref{lem:random}(c) applies to
Theorems~\ref{thm:main} and~\ref{thm:zero-one}.  For $k\ge3$, drawing the
directions as in (i) is not enough for this argument.  A polynomial of
Lemma~\ref{lem:rigidity} can have degree more than one in a single direction,
one for each edge of that colour in the trees, and a distribution can give no
mass to hyperplanes through the origin and yet be carried by a quadric cone, as
the uniform distribution on a small circle of the unit sphere is.  We do not
know whether Theorem~\ref{thm:main} holds for directions drawn in this way.
Rounded coordinates, as in (iv), satisfy none of (i)--(iii); for instance, two
rounded directions are parallel with positive probability.

\paragraph{Which result uses which.}  Table~\ref{tab:uses} summarises, and
Figure~\ref{fig:conditions} shows how the conditions and the lists that
separate fit together.

\begin{table}[ht]
\centering
\small
\begin{tabular}{@{}p{6.6cm}p{8.2cm}@{}}
  \toprule
  Result & Condition on the directions\\
  \midrule
  Theorems~\ref{thm:three} and~\ref{thm:plane}, Lemmas~\ref{lem:pad}%
    --\ref{lem:hypercube}, \ref{lem:criterion}, \ref{lem:four}
    and~\ref{lem:five} & none: they hold
    for every list\\
  Lemma~\ref{lem:linear} (the linear bound) & general position
    (Definition~\ref{def:gp})\\
  Corollary~\ref{cor:generic-no-pattern} and Theorem~\ref{thm:main} & generic
    lists, and so almost every list (Definition~\ref{def:ae})\\
  Theorem~\ref{thm:zero-one} & generic lists, and lists with algebraically
    independent coordinates (Definition~\ref{def:ai})\\
  Corollary~\ref{cor:nested} & random, drawn as in (ii) or (iii) of
    Table~\ref{tab:random}\\
  Corollary~\ref{cor:axes} & the axes fixed, and almost every choice of the
    other directions\\
  Remarks~\ref{rem:small} and~\ref{rem:gp-sharp} & general position, and
    almost every list within a family\\
  Remark~\ref{rem:compress} & a list that separates, and almost every linear
    map\\
  \bottomrule
\end{tabular}
\caption{The conditions that each result puts on the directions.}
\label{tab:uses}
\end{table}

\clearpage
\section*{Supplementary tables}
\setcounter{table}{0}
\renewcommand{\thetable}{S\arabic{table}}
\renewcommand{\theHtable}{S\arabic{table}}

\begin{table}[ht]
\centering
\begin{tabular}{@{}lrrrrr@{}}
  \toprule
  & $n=10^2$ & $10^3$ & $10^4$ & $10^6$ & $10^9$\\
  \midrule
  \multicolumn{6}{@{}l}{$k=2$}\\
  every list fails below & $11$ & $16$ & $20$ & $29$ & $43$\\
  succeed from, before this paper & $101$ & $1001$ & $10\,001$ & $45\,854^*$
    & $147\,481^*$\\
  $m^{\mathrm{HS}}_k(n)$, from \cite{HS26} as written & $364$ & $558$ & $764$
    & $1198$ & $1890$\\
  $m^{\mathrm{HSE}}_k(n)$, equation~\eqref{eq:mhse} & $288$ & $458$ & $638$ & $1024$ & $1646$\\
  $m^{\mathrm{CF}}_k(n)$, equation~\eqref{eq:mcf} & $90$ & $144$ & $203$ & $329$ & $533$\\
  $m^\sharp_k(n)$, equation~\eqref{eq:msharp} & $82$ & $134$ & $190$ & $310$ & $508$\\
  $m^\flat_{k,\mathrm{eq}}(n)$, equation~\eqref{eq:mflateq} & $82$ & $134$ & $190$ & $310$ & $508$\\
  $m^\flat_{k,\mathrm{bal}}(n)$, equation~\eqref{eq:mflatbal} & $82$ & $133$ & $189$ & $310$ & $507$\\
  \addlinespace
  \multicolumn{6}{@{}l}{$k=3$}\\
  every list fails below & $16$ & $22$ & $29$ & $42$ & $62$\\
  succeed from, before this paper & $201$ & $2001$ & $20\,001$
    & $2\cdot10^6{+}1$ & $2\cdot10^9{+}1$\\
  $m^{\mathrm{HS}}_k(n)$, from \cite{HS26} as written & $597$ & $912$ & $1242$
    & $1938$ & $3045$\\
  $m^{\mathrm{HSE}}_k(n)$, equation~\eqref{eq:mhse} & $483$ & $759$ & $1053$ & $1677$ & $2679$\\
  $m^{\mathrm{CF}}_k(n)$, equation~\eqref{eq:mcf} & $170$ & $270$ & $377$ & $605$ & $975$\\
  $m^\sharp_k(n)$, equation~\eqref{eq:msharp} & $162$ & $258$ & $363$ & $588$ & $951$\\
  $m^\flat_{k,\mathrm{eq}}(n)$, equation~\eqref{eq:mflateq} & $162$ & $258$ & $363$ & $588$ & $951$\\
  $m^\flat_{k,\mathrm{bal}}(n)$, equation~\eqref{eq:mflatbal} & $160$ & $258$ & $363$ & $587$ & $951$\\
  \bottomrule
\end{tabular}
\caption{Numbers of directions for separating multisets of $n$ points in $\R^2$
and $\R^3$, as in Figure~\ref{fig:best}.  Every list of fewer directions than the
first row fails: for $k=2$ by Theorem~\ref{thm:plane}; for
$k=3$ by Theorem~\ref{thm:three} and the logarithmic bound (before
Theorem~\ref{thm:three}, one fewer, except at $n=10^4$).  The second row is the
smallest number of directions from which lists were known to succeed before
this paper: the linear bound $(k-1)n+1$, for every list in general position
\cite{Renyi,HSW}, or, for the starred entries, the bound $5\log_2^3(2n)$ of
\cite[Proposition~3.1]{MPS} with the constant that they state, for almost every
list.  The remaining rows are for
almost every list: Hair and Sahai's lemmas as their proof writes them
(Section~\ref{sec:graphs}), $m^{\mathrm{HSE}}_k(n)$, which is their final
inequality evaluated exactly (equation~\eqref{eq:mhse}), and the four bounds of Theorem~\ref{thm:main}.
$m^\flat_{k,\mathrm{bal}}(n)$ is the smallest known upper bound, and $m^{\mathrm{CF}}_k(n)$ the
smallest known in closed form for these $k$.  Printed, with the credits, by
\texttt{boundsdb.py --table} in the accompanying code; \texttt{flatbound.py}
checks the rows for $m^\flat_{k,\mathrm{eq}}(n)$ and
$m^\flat_{k,\mathrm{bal}}(n)$ by a direct search.}
\label{tab:numbers}
\end{table}

\begin{table}[ht]
\centering
\begin{tabular}{@{}lrrrr@{}}
  \toprule
  $n$ & $k$ & $m^\sharp_k(n)$ & $m^\flat_{k,\mathrm{eq}}(n)$ & $m^\flat_{k,\mathrm{bal}}(n)$\\
  \midrule
  $10^2$ & $2$  & $82$   & $82$   & $82$\\
         & $3$  & $162$  & $162$  & $160$\\
         & $8$  & $640$  & $630$  & $622$\\
         & $10$ & $860$  & $814$  & $812$\\
         & $20$ & $2040$ & $1752$ & $1751$\\
  \addlinespace
  $10^3$ & $2$  & $134$  & $134$  & $133$\\
         & $3$  & $258$  & $258$  & $258$\\
         & $8$  & $1008$ & $990$  & $990$\\
         & $10$ & $1340$ & $1287$ & $1286$\\
         & $20$ & $3160$ & $2800$ & $2785$\\
  \addlinespace
  $10^4$ & $2$  & $190$  & $190$  & $189$\\
         & $3$  & $363$  & $363$  & $363$\\
         & $8$  & $1392$ & $1386$ & $1383$\\
         & $10$ & $1850$ & $1793$ & $1791$\\
         & $20$ & $4320$ & $3887$ & $3879$\\
  \addlinespace
  $10^6$ & $2$  & $310$  & $310$  & $310$\\
         & $3$  & $588$  & $588$  & $587$\\
         & $8$  & $2208$ & $2208$ & $2206$\\
         & $10$ & $2920$ & $2871$ & $2865$\\
         & $20$ & $6780$ & $6210$ & $6200$\\
  \bottomrule
\end{tabular}
\caption{Numerical comparison of the bound $m^\sharp_k(n)$
(equation~\eqref{eq:msharp}), in which all $k$ blocks must connect $W$, the equal-block bound $m^\flat_{k,\mathrm{eq}}(n)$
(equation~\eqref{eq:mflateq}), and the balanced-block bound
$m^\flat_{k,\mathrm{bal}}(n)$ (equation~\eqref{eq:mflatbal}) of
Theorem~\ref{thm:main}, for $n\in\{10^2,10^3,10^4,10^6\}$ and
$k\in\{2,3,8,10,20\}$.  Allowing $q>k$ blocks and requiring only $k$ of them to
connect $W$ reduces the number of directions needed, with larger savings as the
dimension $k$ grows.  Computed by \texttt{flatbound.py} in the accompanying code.}
\label{tab:flat-comparison}
\end{table}

\newpage
\section{Index: which results exclude which values of $m_k(n)$ ?}
\label{app:exclusion}

\begin{figure}[!htbp]
\centering
\includegraphics[width=\textwidth]{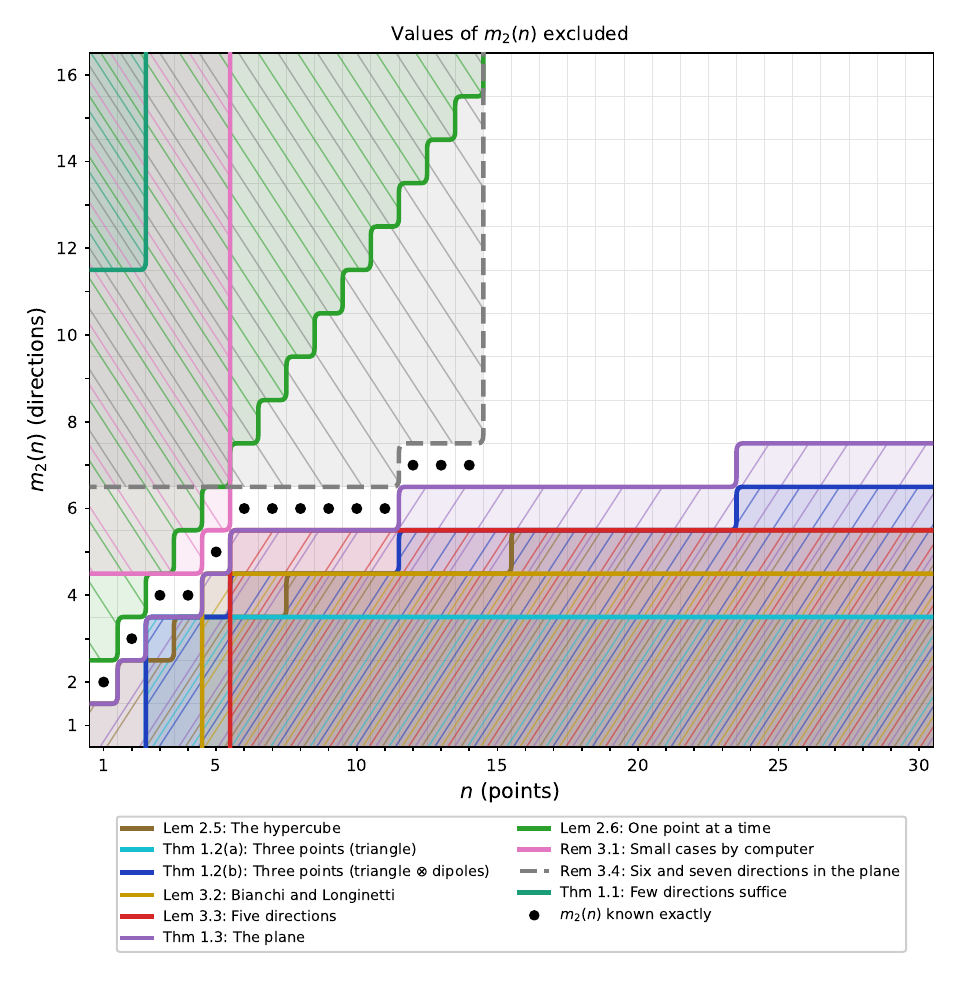}
\caption{All $k=2$ results together.  Dots mark the $n$ at which only one
value of $m_2(n)$ is left.}
\label{fig:excl-k2-combined}
\end{figure}

\begin{figure}[p]
\centering
\includegraphics[width=\textwidth]{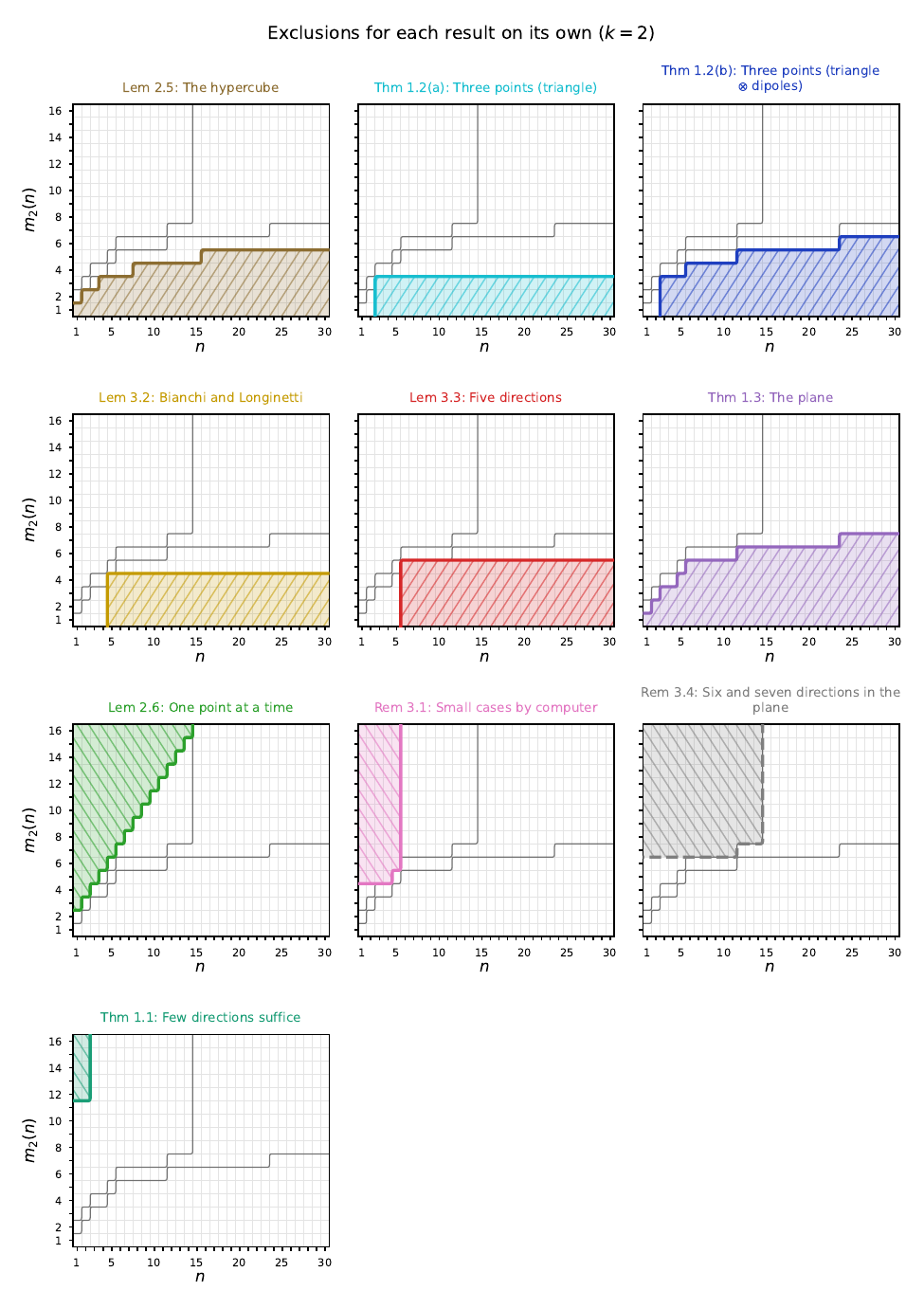}
\caption{Exclusions for each result on its own for $k=2$.  The thin black line
outlines what all of them together leave.}
\label{fig:excl-k2-singles}
\end{figure}

\begin{figure}[p]
\centering
\includegraphics[width=\textwidth]{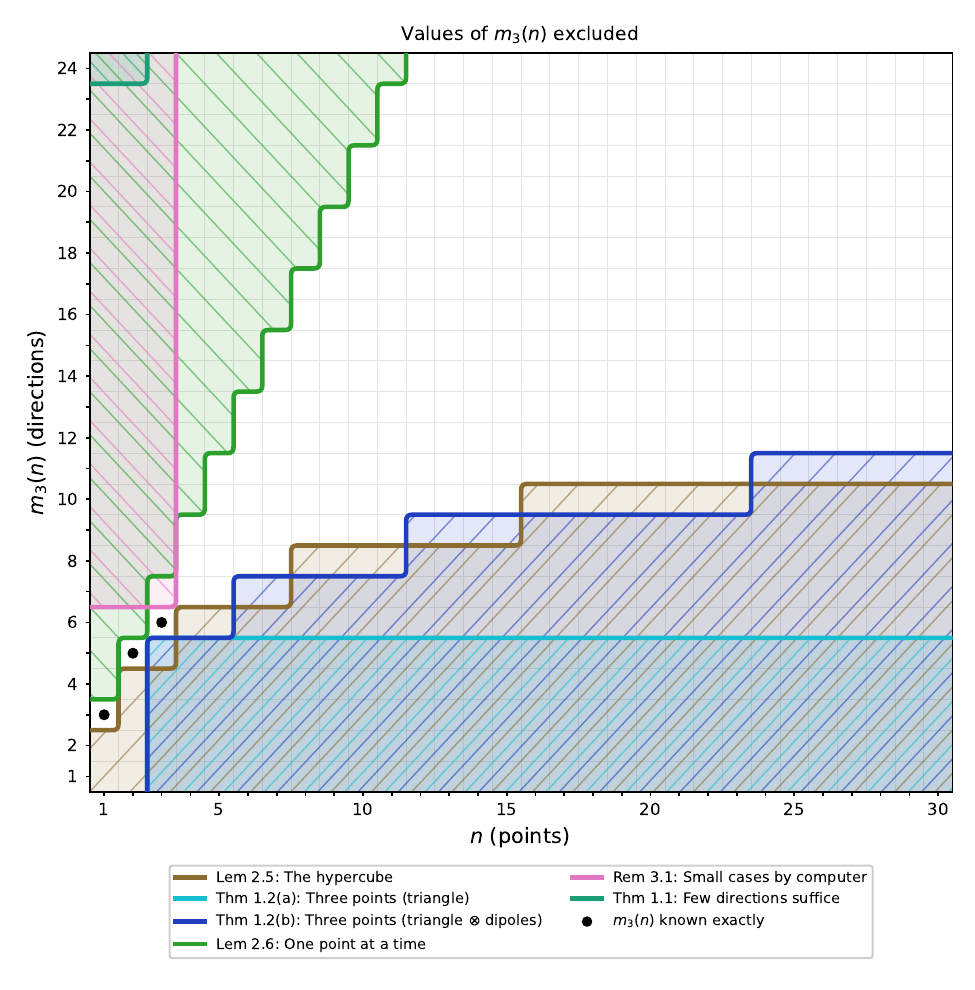}
\caption{All $k=3$ results together.  Dots mark the $n$ at which only one
value of $m_3(n)$ is left.}
\label{fig:excl-k3-combined}
\end{figure}

\begin{figure}[p]
\centering
\includegraphics[width=\textwidth]{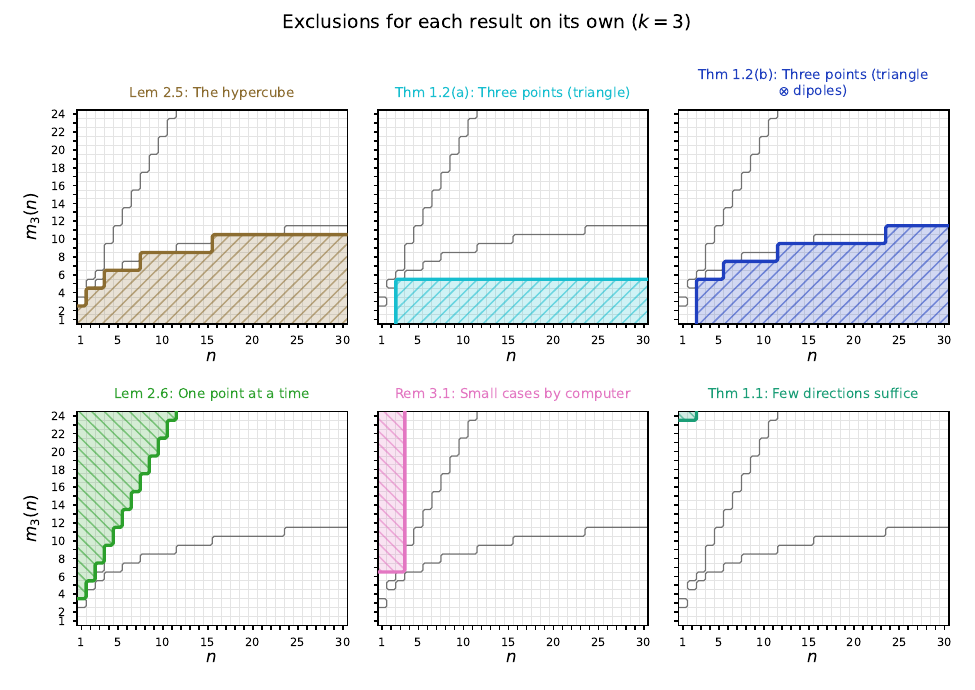}
\caption{Exclusions for each result on its own for $k=3$.  The thin black line
outlines what all of them together leave.}
\label{fig:excl-k3-singles}
\end{figure}

Figures~\ref{fig:excl-k2-combined} to~\ref{fig:excl-k3-singles} show,
for $k=2$ and $k=3$ and $n\le30$,
the values of $m_k(n)$ which are ruled out by each of the main results reported in this note.
This provides a mechanism for comparing the differing natures of the contributions made by each result.
Each cell $(n,m)$ stands for the statement $m_k(n)=m$.  Two facts extend every
result for free, both stated after Lemma~\ref{lem:pad}: $m_k$ is nondecreasing
in $n$, and adding directions never hurts.  So a lower bound $m_k(n_0)\ge L$
excludes every cell with $n\ge n_0$ and $m<L$ (down and to the right,
hatched~/), and an upper bound $m_k(n_0)\le U$ excludes every cell with
$n\le n_0$ and $m>U$ (up and to the left, hatched~\textbackslash).  
The exclusion regions of results that combine other results are drawn only if the new exclusion region is not a simple
union of those contributing; Theorems~\ref{thm:three}(b) and~\ref{thm:plane} qualify.
Results about particular lists, such as Remark~\ref{rem:gp-sharp}, exclude
nothing, since $m_k$ is a minimum over lists.  Remark~\ref{rem:six-seven} is
dashed, as its proofs rest on the correctness of the search code.
Exclusions from $m^\sharp_k(n)$, $m^\flat_{k,\mathrm{bal}}(n)$ and the other bounds
from Theorem~\ref{thm:main} are not visible on most of these plots, as they become
important only at larger $n$.
Each result is labelled by its number
and title or a short description.

\end{document}